\documentclass[12pt]{article}
      
\usepackage{ascmac,amssymb,amsthm, graphicx,latexsym,amsmath}

\newtheorem{theorem}{Theorem}[section]
\newtheorem{lemma}{Lemma}[section]
\newtheorem{remark}{Remark}[section]

\usepackage{color}

\def\green{\color{green}}

\def\0{\boldsymbol{0}}

\def\CC{\mathbb{C}}

\makeatletter

\@addtoreset{equation}{section}

\makeatother

\begin{document}

\begin{center}
{\bf  \Large{Pseudoconvex domains in the Hopf surface}}
\end{center}

\begin{center} by 
\end{center}
\begin{center}
 \large{ Norman Levenberg and Hiroshi Yamaguchi}\ 
\end{center}

\section{Introduction} Let $a\in \mathbb{ C} \setminus \{0\} $ with  $|a|>1$ 
and 
let { ${\mathbb  H}_a$} be the {\sl Hopf manifold with respect to $a$}, i.e., 
${  {\mathbb  H}_a}= \mathbb{ C}^n \setminus \{0\}/ \sim$ where $ z' 
\sim  z$ if and only if there exists $ m\in  \mathbb{ Z}$ such that 
$ z'=a^m  z$ in $\mathbb{ C}^n \setminus \{0\}$.   
In a previous paper \cite{kly} we showed that any pseudoconvex 
domain $D \subset {  {\mathbb  H}_a}$ with $C^\omega-$smooth boundary which 
is not Stein 
is biholomorphic to $T_a \times D_0$ where $D_0 $ is a Stein domain 
in $\mathbb{ P}^{n-1} $ with $C^\omega-$smooth boundary and $T_a$ is a one-dimensional torus. 
This was achieved using the technique of variation of domains in a 
complex Lie group developed in \cite{kly} applied to { ${\mathbb  H}_a$} as a complex homogeneous space with transformation group $GL(n,\mathbb{C})$ (Theorem 6.5 in \cite{kly}).

For $a, b\in \mathbb{ C} \setminus \{0\}$ with $|b|\ge |a|>1$ 
we let ${\mathcal  H}_{(a,b)}$ be the {\sl Hopf surface with respect to $(a,b)$}, i.e., 
${\mathcal  H}_{(a,b)} = \mathbb{ C}^2 \setminus \{(0,0)\}/ \sim$, where $(z,w) 
\sim (z',w')$ if and only if there exists $n\in  \mathbb{ Z}$ such that $z'=a^nz, \ w'=b^nw$. 
We set ${\bf  T}_a=T_a \times \{0\}$, ${\bf  T}_b=\{0\} \times T_b$, and 
${\mathcal  H}^* = {\mathcal  H} \setminus ({\bf  T}_a \cup {\bf  T}_b)$. 
We define 
\begin{equation}\label{rhoeqn}
\rho:= \frac{ \log |b|}{\log |a|}  \ge 1.
\end{equation}
We remark that ${\mathcal  H}_{(a,b)}$ is not a complex Lie group. { 
 However, ${\mathcal  H}^*$ is both a complex Lie group and a complex 
homogeneous space.} With the aid of the aforementioned technique of variation of domains in \cite{kly}, we can characterize the pseudoconvex domains with $C^\omega-$smooth boundary in ${\mathcal  H}_{(a,b)}$ 
 which are not Stein. 
\begin{theorem} \label{thm-1}  
  Let $D$ be a pseudoconvex 
domain in ${\mathcal  H}_{(a,b)}$ with $C^\omega-$smooth boundary. Suppose $D$ is not Stein.
\begin{enumerate}
\item [{\rm Case {\bf  a}:}] \   $\rho $ is irrational. 

 We set 
       $\Sigma _c=\{ |w|= c|z|^\rho\}/ \sim$ for $c\in (0, +\infty)$, 
      where $\{ |w|= c|z|^\rho\}= \{|w|=c|z|^\rho: z \in \mathbb{ C}^*\} 
      \subset \mathbb{ C}^* \times \mathbb{ C}^*\}$,  
$\Sigma _0= {\bf  T}_a$ and $\Sigma _{+\infty}={\bf  T}_b$. We have  
${ {\mathcal  H}= \bigcup_{c\in [0, +\infty] }\Sigma _c}$ where this is a 
disjoint union and each 
       $\Sigma _c$ for $c\in (0, \infty)$ is a compact Levi flat hypersurface in 
       ${\mathcal  H}^*$ biholomorphic to $\Sigma _1$ in 
      ${\mathcal  H}$. Then $D$ reduces 
      to one of the following:
 \begin{enumerate}
 \item [{\rm (a-1)}] \ \  There exist $0<k_1<k_2< +\infty$ such 
       that $
{ D= \bigcup_{c\in (k_1, k_2)}\ \Sigma_c.}
$
 \item [{\rm(a-2')}]  \ \  There exists a positive number $k$ 
such that
$ { 
D=\bigcup_{c\in [0, k)} \ \Sigma _c.}
$
 \item [{\rm(a-2'')}] \ \ There exists a positive number $k$ 
such that
$
{ D=\bigcup_{c\in (k, +\infty]} \ \Sigma_c}.
$
\end{enumerate}
 \item [{\rm Case {\bf b}:}] \  $\rho=q/p$ is rational where 
       $q>p\ge 1$ and $(p,q)=1$. 
Setting 
\begin{equation}\label{taueqn}
\tau: = \frac{ 1}{2\pi}\, (\frac{ q}p \arg\, a - \arg\, b), \quad 0\le 
       \arg\,a, \ \arg\,b <2\pi,
\end{equation}

we have two cases: 
\begin{enumerate}
 \item [{\rm Case {(b1)}:}] \  $\tau$ is irrational. We have the 
       same disjoint union ${ {\mathcal  H}= \bigcup_{[0, +\infty] } 
       \Sigma_c}$  as in Case {\bf  a}, and the domain 
$D$ reduces to 
  one of {\rm(a-1)}, {\rm(a-2')} or {\rm(a-2'')}.

\vspace{2mm} 
 
\item [{\rm Case {(b2):}}] \  $\tau =m/l$ is rational with $l\ge 1$ and $(l,m)=\pm 
       1$ or $\tau=0$ (and we set $l=1$). {  Let $g$ be the greatest common divisor of 
       $p$ and $l$, and set $\nu:= pl/g \in \mathbb{ Z}$.}
We define $K:= {  \{ e^{2\pi i\,  k / \nu }\}_{k=0, 1,\ldots 
       , \nu-1}}$, a subgroup of $\mathbb{ C}^*$, and define 
       $\sigma_c:= \{w= c\, z^\rho\}/ \sim$ for $c\in \mathbb{ C}^*/K$, 
       where 
       $\{w= c\, z^\rho\}=\{ w= c z^\rho: z\in  \mathbb{ C}^*\} \subset \mathbb{ C}^* \times 
       \mathbb{ C}^*$. Setting $\sigma _0= {\bf  T}_a$ and 
$\sigma _{\infty}={\bf  T}_b$, we have ${\mathcal  H} = \bigcup_{c\in 
       \mathbb{ P}^1}  \sigma_c$ where this is a disjoint union. There exists a domain $\delta $ in $\mathbb{  P}^1$ with 
       smooth boundary such that $ { 
D=\bigcup_{c\in \delta }\, \sigma_c.}$

\end{enumerate}\end{enumerate}\end{theorem}

\noindent In Case (b2), $\sigma_1 
= \{ w= z^\rho\}/ \sim$ is a compact curve in ${\mathcal  H}^*$ which is 
       a subgroup of ${\mathcal  H}^* $ and which, as a Riemann surface, is equivalent to 
a torus. Each $\sigma _c$ for $c\in \mathbb{ C}^*$ is 
biholomorphic to $\sigma _1$ in $ {\mathcal  H}^*$. We consider $\mathbb{ C}^* /K$ as a Riemann surface which is equivalent to $\mathbb{ C}^* $.  Then $\{w= c \, z^\rho\}/ \sim\ =\ 
\{w= c' \, z^\rho\}/ \sim$ for $c'/c\in K$. Here we have used the shorthand notation $c= cK\in \mathbb{ C}^*/K$.

 \vspace{2mm} The main idea behind the proof is this: starting with a pseudoconvex domain $D\subset \mathcal H$ with smooth boundary, we consider $D^*=D\cap \mathcal H^*$. We construct a natural plurisubharmonic exhaustion function using our $c-$Robin function techniques in \cite{kly}. It is natural to try to extend this function to $D$ first as a plurisubharmonic function and then as an exhaustion function. We study obstructions to the resulting function (or a modification of it) being strictly plurisubharmonic arising from the possible existence of certain holomorphic vector fields. As a by-product of this procedure, we also encounter an interesting class of Stein subdomains in $\mathcal H$ which we call {\it Nemirovskii-type} domains. 

\vspace{2mm} The outline of our paper is the following. In the next section, we briefly discuss properties of the Hopf surface ${\mathcal  H}_{(a,b)}$, and in section 3 we state without proof some preliminary results, including a classification of the holomorphic vector fields on ${\mathcal  H}_{(a,b)}$ and their integral curves. 
We also indicate why the domains listed in Theorem \ref{thm-1} are not 
Stein. The proof of Theorem \ref{thm-1} is given in section 4. At the end of that section we give an example of the aforementioned Nemirovskii-type domains. The proofs of the results in section 3 are given at the end of the paper in Appendix A and Appendix B. 

\vspace{2mm} We would like to thank Professor Tetsuo Ueda for suggesting this problem. 
\section {Properties of the Hopf surface ${\mathcal  H}_{(a,b)}$}
We write $\CC^* :=\CC \setminus \{0\}$ and $(\CC^2)^*:=
 \CC^2 \setminus \{(0,0)\}$.
Fix $a,b \in \CC^*$ with $1<|a|\le|b|$. 
For  $(z,w), 
 (z',w') \in (\CC^2)^*$, we define the equivalence relation 
$$
(z,w) \sim (z',w') \quad \mbox{ iff } \quad \exists \ n \in  \mathbb{ Z}
 \mbox{ such that  } \ z'= a^n z, \ w'= b^n w.
$$
The space  $(\CC^2)^* / \sim $ consisting of all equivalence classes
$$
[z,w]: = \{(a^nz, b^nw): n\in  \mathbb{ Z} \}, \quad (z,w)\in (\CC^2)^*
$$
is called  
the {\it  Hopf surface}   
${\mathcal  H}={\mathcal H}_{(a,b)}$; it is a complex two-dimensional 
compact manifold. 

For $z,\ z'\in \mathbb{ C}^*  $ we define $z \sim _a z'$ if and only if there exists $ n\in \mathbb{ Z}$ such that  $z'=a^nz$ in $\mathbb{ C} ^*$. Then 
$$
T_a :=\mathbb{ C}^*  /\sim_a \ \hbox{and} \ T_b :=\mathbb{ C}^*  /\sim_b
$$
are complex one-dimensional tori, and ${\mathcal  H}$ contains two disjoint compact 
analytic curves ${\bf  T}_a=T_a \times 
\{0\}$ and ${\bf  T}_b=\{0\}\times T_b$. We have ${\bf  T}_a \cup {\bf  T}_b= \{(z,w)\in (\mathbb{ C}^2)^*   : zw=0\}/ \sim 
$
 in ${\mathcal  H}$; for simplicity we write ${\mathbf T}_a \cup {\mathbf  T} 
_b = \{zw=0\}$. 
We consider the subdomain ${\mathcal  H}^*$ of ${\mathcal  H}$ 
 defined by 
\begin{equation}\label{hstar}
{\mathcal  H}^* := {\mathcal  H} \setminus \{zw=0\}.
\end{equation}
Thus ${\mathcal  H}$ is a compactification of ${\mathcal  H}^*$ by two 
disjoint one-dimensional tori. The set ${\mathcal  H}^*$ is a complex Lie group 
and will play a crucial role in this work.
\vspace{3mm}

We give a more precise description of the Hopf surface. 
A fundamental domain for ${\mathcal  H}$ is 
\begin{align}
\label{eqn:found}  \nonumber 
 {\mathcal F}:&= (\{|z|\le |a|\}\times
\{|w|\le |b| \} )\setminus ( \{|z|\le1\}\times
\{|w|\le 1 \}) \ \ \\  
&= E_1 \cup E_2 \Subset (\mathbb{ C}^2)^*, 
\end{align}
where 
\begin{alignat*}{3}
 E_1&=E_1 ' \times E_1'':= \{|z|\le |a|\} \times \{1 {  <} |w| \le |b|\},\\
E_2&=E_2' \times E_2'':=
 \{1 {  <} |z|\le |a|\} \times \{ |w|\le |b|\}.
\end{alignat*}
For $k=0,\pm 1, \ldots $ we set ${\mathcal  F}_k:= {\mathcal  F}\ 
\times (a^k,b^k)$. Then ${\mathcal  F}_0= {\mathcal  F}$; each 
${\mathcal  F}_k$ is a fundamental domain; and we have the disjoint union $(\mathbb{ 
C}^2)^*= \cup_{n=-\infty}^{\infty} {\mathcal  F}_n$.

\vspace{2mm} 
The Hopf surface ${\mathcal  H}$ is obtained by gluing the boundaries of $\partial 
{\mathcal  F}$ in the following way:
setting 
\begin{alignat*}{3}
L'_a&:= \{|z|\le |a|\} \times \{|w|=|b|\}, \quad L' _1=\{|z|\le 1\} \times \{|w|=1\};\\
  L''_b&:= \{|z|=|a|\} \times \{|w|\le |b|\}, \quad L''_1=\{|z|= 1\} 
 \times \{|w|\le 1\},
\end{alignat*}
we have the identifications:
\begin{align*}
(1)\ \  &(z, w)\in L_a'  \quad \mbox{ with} \quad (z/a, w/b)\in L'_1;\\
(2)\ \ &(z,w)\in L_b''  \quad \mbox{ with} \quad 
 (z/a, w/b)\in L''_1.
\end{align*}


We set
\begin{eqnarray}
 \label{ai}
{\mathcal  I}=\{(a^n,b^n) \in \mathbb{ C}^* \times  \mathbb{ C}^* : n\in 
 \mathbb{ Z}\} \subset \mathbb{ C}^* \times  \mathbb{ C}^* ,
\end{eqnarray}
which is a discrete  set in $\mathbb{ C}^* \times  \mathbb{ C}^* $.

For a set $D \subset {\mathcal  H}$ we will often simply describe $D$ as a set of points 
in $(\mathbb{ C}^2)^*$ where the equivalence relation $\sim$ is understood. If there is possibility of confusion we will write
\begin{align}\label{tilded}
 \widetilde D &=\{(z,w)\in (\mathbb{ C}^2)^*: [z,w]\in D\} \subset (\mathbb{ 
C}^2)^*, 
\end{align}
so that $\widetilde D= \widetilde D \times {\mathcal  I}$ and hence $\widetilde D/ \sim \ = D$.


We give an example of the action of the equivalence relation  
which will illustrate the difference between the Lie group ${\mathcal  
H}^*$ and the Hopf surface ${\mathcal  
H}$. Let $D=\mathbb{ C}_z \times \{w\}$ where $w\ne 0$. As a subset of ${\mathcal  
H}^*$, the complex curve $D \cap (\mathbb{C}^* \times \mathbb{C}^*) /\sim$ is not relatively compact and is equivalent to $\mathbb{ C}^*$. 
However, as a complex curve in ${\mathcal  
H}$, $D/\sim$ is not closed and is equivalent to $\mathbb{ C}$. Moreover, 
if $|b|^{k-1}<|w|< |b|^k$, then $(0,w)\in {\mathcal  F}_k$ and   
$$
D /\sim \ = D_0 \cup D_1 \cup D_2 \cup \cdots 
$$
where 
$$
D_0= \{|z|< |a|^k\} \times \{w\}, \ \ D_n=\{|a|^{k-1} \le |z|\le |a|^k\} \times \{{ 
w}/{b^{n}} \}, \quad n=1,2,\ldots
$$
Thus $D_0$ is a disk and $D_n, \ n=1,2,\ldots $ are annuli such that 
 $D_{n+1}=D_n \times (1,1/b), \ n=1,2,\ldots $. Hence the 
$D_n, \ n=1,2,3,\ldots $ are conformally equivalent and, as $n\to \infty$, they wind 
around and converge to ${\bf  T}_a$ in ${\mathcal  H}$.

\vspace{3mm} 
Following T.\,Ueda, we consider the following real-valued function $U[z,w]$ on ${\mathcal  
H}^*$: 
$$
U[z,w]= \frac{ \log |z|}{\log |a|} - \frac{\log  |w|}{\log |b|} \quad 
\mbox{ for } \ [z,w]\in {\mathcal  H}^*.
$$
This has the following properties:
\begin{enumerate}
 \item [(1)] $U[z,w]$ is a pluriharmonic function on ${\mathcal  H}^*$ 
       satisfying  
$$
\lim_{ [z,w]\to {\bf T}_a} U[z,w]=+\infty \quad \mbox{ and } \quad  
\lim_{[z,w]\to {\bf  T}_b} U[z,w]=-\infty,
$$
thus for any interval $I \Subset (-\infty, \infty)$,  
the subdomain $U^{-1}(I)$ of ${\mathcal  H}^*$ is relatively compact in ${\mathcal  H}^*$.
 \item [(2)]  $|U[z,w]|:={\rm Max} \{U[z,w], -U[z,w]\} $ is a plurisubharmonic 
       exhaustion function for ${\mathcal  H}^*$ which is pluriharmonic 
       everywhere except on the Levi-flat set 
$$
\frac{ \log |z|}{\log |a|} = \frac{ \log |w|}{\log |b|}, \quad 
       i.e., \quad |w|=|z|^\rho \ \   \mbox{ in } \ \ {\mathcal  H}^*.
$$
 \item [(3)] For $c \in (-\infty, +\infty)$, 
 the level set 
$$
S_c: \quad U[z,w]=c
$$
is equal to $|w|= k |z|^\rho$ where $k= e^{-c \log |b|}>0$.  Thus $\{ k_1|z|^\rho\le |w|\le k_2|z|^\rho\}$ is equal to 
$U^{-1}([c_1,c_2])$ where $k_i=e^{-c_i\log \ |b|}$; 
    \ 
while $\{|w|\le k|z|^\rho\}$ is equal to 
$U^{-1}([c,+\infty)) \cup {\bf  T}_a$; and 
$\{|w|\ge k|z|^\rho\}$ is equal to 
$U^{-1}((-\infty, c])) \cup {\bf  T}_b$ where $k=e^{-c\log \, |b|}$.
\end{enumerate}

From (2) and (3), each of the domains $D$ in ({\rm a}-1), ({\rm a-2'})
and ({\rm a}-2'') in the statement of Theorem \ref{thm-1} 
contains a compact, Levi-flat hypersurface $S_c$ for appropriate $c$; 
hence each such $D$ is not Stein. 

\section{Preliminary results}

In this section, we discuss two basic results which we will need. The first concerns holomorphic vector fields in ${\mathcal  H}$, while the second concerns general pseudoconvex domains with $C^\omega-$smooth boundary in $\mathbb{ C}^2$. 

We consider the linear space $\mathfrak{ X}$ of all holomorphic vector fields $X$ of the form 
$$
X= \alpha z \frac{\partial }{\partial z} + \beta w \frac{\partial 
}{\partial w}, \quad  \alpha , \beta \in \mathbb{ C}
$$
in $\mathbb{ C}^* \times  \mathbb{ C}^* $. Any such $X$ clearly induces a holomorphic vector field on 
${\mathcal  H}$. The integral curve of $X$ with 
initial value $(z_0, w_0) \in {  \mathbb{ C}^* \times \mathbb{ C}^*}$ is 
$$
(z_0,w_0) \exp tX= \left\{ \ \begin{array}{lll}
 z= z_0 e^{\alpha t},\\[2mm]
w=w_0 e^{\beta t}, \end{array} \right. 
\qquad t\in \mathbb{ C}.
$$
Therefore, if, for example, $\alpha \ne 0$, we can write
$$
w= c_0 z^{\beta/\alpha}\quad \mbox{ where $c_0=w_0\bigl/ 
z_0^{\beta/\alpha}$.}
$$
The integral curve 
$[z_0, w_0] \exp tX $ of $X$ in ${\mathcal  H}$ with initial value 
$[z_0,w_0]$ is equal to $ \{w=c_0\, z^{\beta/\alpha}\}/ \sim $ in ${\mathcal  
H}^*$.

In particular, we consider
\begin{equation} \label{xu}
 X_{u}: = (\log|a|)\, z\, \frac{\partial }{\partial z} + (\log|b|) \,
w\,\frac{\partial }{\partial w}.\end{equation}
The integral curve of $X_{u}$ with initial value $(1,1)$ is
\begin{align}\exp {  t} X_u&=  \left\{ \ \begin{array}{lll}
 z=  e^{(\log|a|) t}, \nonumber \\[1mm]
w= e^{(\log|b|) t}{ ,}
 \end{array} \right.
\qquad t\in \mathbb{ C}.
\end{align}
Thus $w= z^\rho$.
We set $\widetilde \sigma _u:=   \{\exp tX_u: t\in \mathbb{ C}\} { / \sim} \subset 
       {\mathcal  H}^* $ and denote by $\tilde \Sigma _u$  the closure of $\widetilde \sigma _u$ 
       in ${\mathcal  H}$. For future use, we define the linear subspace $\mathfrak{ X}_u=\{cX_u: \, c\in \mathbb{ C}\}$ of $\mathfrak{ X}$.

The next lemma gives more precise information about the integral curves and 
will be crucial in the proof of Lemma \ref{essen-lemm}. Recall for rational $\rho= \frac{ \log 
|b|}{\log |a|}$, we write $\rho=q/p, \ p\ge 1, \ (p,q)=1$ and  
 $\tau:= ((q/p) \arg a- \arg b)/2\pi$ was defined in (\ref{taueqn}). 

\begin{lemma} \label{fund-1} 
\begin{enumerate}
\item[1.] For $X_u= (\log|a|)\, z\, \frac{\partial }{\partial z} + (\log|b|) \,
w\,\frac{\partial }{\partial w}$ we have:
\begin{enumerate}
 \item [(1)] In case $\rho $ is irrational or $\tau$ is irrational,  
       $\tilde \Sigma _u=\{|w|=|z|^\rho\} { / \sim }$ is 
a real three-dimensional Levi-flat closed hypersurface
in ${\mathcal  H}^*$ with $\tilde \Sigma _u \Subset {\mathcal  H}^*$.
 \item [(2)] If $\tau$ is rational, then $\widetilde \sigma _u $ is a 
       compact curve in ${\mathcal  H}^*$ which, as a Riemann surface, is equivalent to a  
       one-dimensional torus. 
\end{enumerate}
\item [2.] For $X=\alpha z \frac{\partial }{\partial z} + \beta w \frac{\partial 
}{\partial w} \not\in \{cX_u: c\in \mathbb{ C}\}$, the integral 
	  curve  $\sigma :=\{\exp \,tX: \, t\in \mathbb{ C}\}{ / \sim}$ in 
	  ${\mathcal  H}^*$ is 
	  not relatively compact in ${\mathcal  H}^*$. If we let $\Sigma $ denote the closure of 
$\sigma $ in ${\mathcal  H}$, then: 
\begin{enumerate}
 \item [(1)] If $\alpha, \beta \ne 0$, we have  
$\Sigma  \supset {\bf  T}_a \cup {\bf  T}_b$.
 \item [(2)] If only one of $\alpha$ or $\beta$ is not $0$, 
       e.g., $\alpha \ne  0$ and $\beta= 0$, we have 
$\Sigma \supset {\bf T}_a$ and $\Sigma \cap {\bf T}_b=\emptyset$.
\end{enumerate}\end{enumerate}\end{lemma}

\noindent The proof of Lemma \ref{fund-1} is in Appendix A. 

{ The following conclusions from Lemma \ref{fund-1} and the argument in Appendix A will be needed in proving the key Lemma \ref{essen-lemm}:} 
\begin{enumerate}
 \item [($\alpha$)] If $\rho$ in (\ref{rhoeqn}) is in Case ${\bf  a}$ or 
       Case $({\rm b1})$ in Theorem 
       \ref{thm-1}, i.e., either $\rho$ is irrational or $\rho$ is rational and $\tau$ in (\ref{taueqn}) is irrational, then 
\begin{align*}
 {\mathcal  H}= \bigl(\bigcup_{c\in (0,\infty)} \{|w|=c |z|^\rho\}/\sim \ \bigr) 
\cup ({\bf  T}_a \cup {\bf  T}_b)
\end{align*}
 and this is a disjoint union. Here $\Sigma _c:=\{|w|=c |z|^\rho\}{  /\sim} $ is  
 the closure of the integral curve 
$\sigma[z_0,w_0]= [z_0,w_0] \exp tX_u $ with $c= | w_0/z_0^{\rho}|$ 
       in ${\mathcal  H}$; $\Sigma _c$ is a real three-dimensional Levi flat hypersurface in 
       ${\mathcal H}$ (and hence $\Sigma _c \Subset {\mathcal  H}^*$) 
       which is biholomorphic to $\Sigma _1$ in ${\mathcal  H}$. 
       We set $\Sigma 
       _0 ={\bf  T}_a$ and $\Sigma_\infty={\bf  T}_b$ so that ${\mathcal  H}= 
\bigcup_{c\in [0, \infty]}\,\Sigma_c$.

 \item [($\beta$)] 
If $\rho$ in (\ref{rhoeqn}) is in Case ({\rm b2}) in Theorem 
       \ref{thm-1} so that $\tau$ in (\ref{taueqn}) is rational, 
then for  $[z_0,w_0]\in {\mathcal  H}^*$, the integral curve $[z_0,w_0] 
       \exp tX_u$ in ${\mathcal  H}^*$ is given by $\sigma _c:= \{w=c 
       z^\rho\}/ \sim $ with $ c= w_0/z_0^\rho$ which, we recall, is 
       biholomorphic to the torus $\sigma _1$. Setting $K:=  \{e^{i2\pi 
        k/\nu}\}_{k=0, 1,\ldots ,  \nu-1}$ where $\nu$ was defined in Case ({\rm b2}), we have $\sigma _c= \sigma 
       _{c'}$ if and only if $c/c'\in K$. Considering $\mathbb{ C}^*/K$  
        as Riemann surface equivalent to $\mathbb{ C}^*$ and writing $c=cK$, we 
       have ${\mathcal  H}^*= \bigcup _{c\in \mathbb{ C}^*}\sigma _c$ and this is a disjoint union. 
We note that ${\bf  T}_a= [z_0,0] \exp tX_u$ where $z_0\ne 0$ and  ${\bf  
       T}_b=[0,w_0]\exp tX_u$ where $w_0\ne 0$. 
We set $\sigma 
       _0 ={\bf  T}_a$ and $\sigma_\infty={\bf  T}_b$ so that 
       ${\mathcal  H}=\bigcup_{c\in \mathbb{ P}^1}\, \sigma _c$. 
\end{enumerate}


We now turn to an elementary 
property of a pseudoconvex domain $D$ with 
$C^\omega-$smooth boundary in $\mathbb{ C}^2$. In $\mathbb{ C}^2=\mathbb{ C}_z \times \mathbb{ C}_w$ we consider disks $$
\Delta _1=\{|z|<r_1\}, \quad \Delta _2=\{|w|< r_2\} \quad 
$$
and the bidisk $\Delta =\Delta 
_1 \times \Delta _2$. Let $D$ be a pseudoconvex domain 
 with $C^\omega$ boundary in $\Delta $. We do not assume $D$ is relatively compact. 
Thus there exists a $C^\omega-$smooth, real-valued function 
$\psi(z,w)$ on $\overline{ \Delta }$ such that  
\begin{align*}
 D&= \{(z,w)\in \Delta : \psi(z,w)<0\};\\
\partial D \cap \Delta&= \{(z,w)\in \Delta : \psi(z,w)=0\},
\end{align*}
and on 
$\psi(z,w)=0$ we have both
$\nabla _{(z,w)}\psi(z,w)\ne 0$ and the Levi form ${\mathcal  L}\psi(z,w) \ge 0$. 
We write out this last condition: for 
\begin{align} \nonumber 
  {\mathcal  L}  \psi(z, w)=
 \frac{\partial^2 \psi }
      {\partial z \partial \overline{z}} 
 |\frac{\partial \psi }{\partial w}|^2
-2 \Re\, \left\{
\frac{\partial^2 \psi }{\partial z \partial \overline{ w}}
\frac{\partial \psi}{\partial \overline{ z }} \frac{\partial \psi}{\partial w}
\right\}
+  \frac{\partial^2 \psi }{\partial w \partial \overline{ 
w}}|\frac{\partial \psi }{\partial z}|^2,
\end{align}
\begin{align}\hbox{we have}  \qquad {\mathcal  L} \psi(z, w)\ge 0  \qquad \mbox{on $\psi(z,w)=0$.}\label{levi}
\end{align}
We may assume
$$
\psi(0,0)=0
 \quad \mbox{ and} \quad \frac{\partial \psi }{\partial w}(0,0)\ne 0 
$$
so that $\{w:\psi(0,w)=0\}$ is a $C^\omega-$smooth simple arc in $\Delta_2 $ passing 
through $w=0$.

We set 
$
{\mathcal  S}:=\partial D \cap \Delta,
$
\begin{align*}
 D(z)&:= \{ w\in \Delta _2: (z,w)\in D\} \subset \Delta_2; \ \hbox{and}\\
 S(z)&:= \{w\in \Delta_2: (z,w) \in {\mathcal  S}\} \subset \Delta_2,
\end{align*}
so that $D=\cup_{z\in \Delta _1}(z, D(z)) \subset \Delta $ and ${\mathcal  S}=\cup _{z\in 
 \Delta_1}(z, S(z)) \subset \Delta $. 
Taking $r_1,r_2>0$ sufficiently small we can insure that 
\begin{enumerate}
 \item [(i)] for each $z\in \Delta _1$, $D(z)$ is a non-empty domain 
in $\Delta _2$ and $S(z)$ is a $C^\omega-$smooth open arc in $\Delta 
       _2$ connecting two points $a(z)$ and $b(z)$ on $\partial 
       \Delta _2$; 
 \item [(ii)] $0\in S(0)$.
\end{enumerate}
We also need to assume the following condition for  Lemma 
\ref{ps-lemma}:
\begin{enumerate}
 \item [${(iii)}$] \ \ $\psi (z,0)\not \equiv 0$ in $\Delta 
       _1$, hence, for any disk $\delta_1 =\{|z|<r\} \subset \Delta _1$, there exists $z_0\in \delta_1$  
with $0\not \in S(z_0)$.
\end{enumerate}

 
Under these three conditions we have the following.
 \begin{lemma} \label{ps-lemma} \ \ 
  For any disk $\delta_1=\{|z|<r\}\subset \Delta_1$, there exists a disk $\delta_2=\{|w|<r'\}\subset \Delta_2$ with  
$$
\bigcup _{z\in \delta _1} S(z) \supset D(0)\cap \delta _2.
$$
\end{lemma}

\noindent The proof of Lemma \ref{ps-lemma} is in Appendix B. This result will be used in proving Lemma \ref{lem:basic}.

\section{Construction of the plurisubharmonic exhaustion function 
 $-\lambda [z,w]$ on $D$ }

Let $(\alpha , \beta ) \in \mathbb{ C}^* \times \mathbb{ C}^*$. If we 
 define 
$$
(\alpha ,\beta ): \ [z,w]\in {\mathcal  H} \mapsto 
[\alpha z, 
\beta w] \in {\mathcal  H},
$$
then $(\alpha , \beta )$ is an automorphism of ${\mathcal  H}$. Thus 
$\mathbb{ C}^* \times \mathbb{ C}^*$ acts as a commutative group of automorphisms 
of ${\mathcal  H}$ with identity element $e=(1,1)$. 
Although $\mathbb{ C}^* \times \mathbb{ C}^*  $ is not 
transitive on ${\mathcal  H}$, it is transitive on ${\mathcal  H}^*$. Hence {\it  ${\mathcal  H}^*$ is a complex homogeneous space with Lie transformation group 
 $\mathbb{ C}^* \times  \mathbb{ C}^*$ which acts transitively}. This is the setting of Chapter 6 of \cite{kly}. For any  $[z,w]\in {\mathcal  H}^*$ the isotropy subgroup 
$I_{[z,w]}$ of $ \mathbb{ C}^* \times  \mathbb{ C}^* $ is 
\begin{align*}
 I_{[z,w]}: &=\{(\alpha , \beta )\in \mathbb{ C}^* \times \mathbb{ C}^* : 
(\alpha ,\beta )[z,w]=[z,w]\}\\
&= \{(a^n, b^n)\in \mathbb{ C}^* 
\times  \mathbb{ C}^*: n\in \mathbb{ Z}\}\\
&={\mathcal  I}\mbox{ \  in (\ref{ai})},
\end{align*}
and thus is independent of $[z,w] \in 
{\mathcal  H}^*$. 
We have 
$${\mathcal  H}^*= (\mathbb{C}^* \times \mathbb{C}^* )/ {\mathcal  I} .
$$

In what follows we will generally consider the restriction to $\mathbb{ C}^* \times \mathbb{ C}^*$ 
of the Euclidean metric $ds^2= 
|dz|^2+|dw|^2$ on $\mathbb{ C}^2$, and we fix a positive real-valued 
function $c(z,w)$ of class $C^\omega $ on $\mathbb{ C}^2$. This allows us to define {\it $c-$harmonic functions} and thus a {\it $c-$Green function and $c-$Robin constant} associated to a smoothly bounded domain $\Omega \Subset \mathbb{ C}^* \times \mathbb{ C}^*$ and a point $p_0\in \Omega$ (if $\Omega \not \Subset \mathbb{ C}^* \times \mathbb{ C}^*$ we define these by exhaustion); cf., chapter 1 of \cite{kly}. Varying the point $p_0$ yields the {\it $c-$Robin function} for $\Omega$. However, we remark that any K\"ahler metric $dS^2$ and positive function $C(z,w)$ of class $C^\omega $ on $\mathbb{ C}^* \times \mathbb{ C}^*$ gives rise to a $C-$Green function and hence a $C-$Robin function on $\Omega$; this flexibility will be used in the $4^{th}$ case of the proof of Lemma  \ref{last-lem}. For simplicity, we will always take $c(z,w)$ (or $C(z,w)$) to be a positive constant.

\vspace{3mm} 
{\sl In this section we always assume that $D \subset {\mathcal  H}$ is 
a pseudoconvex domain with $C^\omega-$smooth  
boundary in ${\mathcal  H}$.} Our first goal is to construct a plurisubharmonic exhaustion function for $D$. We note, as observed at the end of section 2, that 
\begin{equation} \nonumber 
\hbox{if} \ D \supset {\bf T}_a \ \hbox{or} \ D \supset {\bf T}_b, 
\ \hbox{then} \ D \ \hbox{is not Stein}. \end{equation}
We define 
$$D^*:=D \cap \{zw\not =0\}\subset {\mathcal  H}^*$$
(see (\ref{hstar})). The distinction between $D\subset {\mathcal  H}$ and $D^*\subset {\mathcal  H}^*$ will be very important. Since $(z, w)\in \mathbb{ C}^*  \times \mathbb{ C}^*  $ defines an 
automorphism of ${\mathcal  H}$,  for $[z,w]\in D$ we can define 
$$
D[z,w]= \{(\alpha ,\beta )\in \mathbb{ C}^* \times \mathbb{ C}^*: 
(\alpha, \beta )[z,w] \in D\}\subset \mathbb{ C}^* \times \mathbb{ C}^*.
$$

Equivalently, using the notation $D \cap {\bf  T}_a= D_a \times 
\{0\}$, $D \cap {\bf  T}_b= \{0\} \times D_b$, 
$\widetilde { D_a}=\{ a^nz: z\in D_a, \ n\in \mathbb{Z}\} \subset \mathbb{ C}_z^*$ and 
$\widetilde {D_b}=\{ b^nw : w\in D_b, \ n\in \mathbb{Z}\} \subset \mathbb{ C}_w^*$, we have

\begin{equation*}
 \begin{array}{lll}
  & D[z,w]= \bigl(\displaystyle{ (\frac{ 1}z, \frac{ 1}w)\cdot D}\bigr) \times {\mathcal 
 I} & \ \  \mbox{if } \ [z,w] \in D^*;\\[2mm]
&D[z,0]= \bigl(\displaystyle{ \frac{ 1}z D_a, \mathbb{ C}^*}\bigr) \times {\mathcal  I}
  \ =(\displaystyle{ \frac{ 1}z}\, \widetilde { D_a}) \times \mathbb{ C}_w^* &  \ \ \mbox{if } \ [z,0] \in  D\cap {\bf T}_a;\\ [2mm]
&D[0,w]= \bigl(\displaystyle{ \mathbb{ C}^*, \frac{ 1}w D_b} \bigr) \times {\mathcal  I}
 \ = \mathbb{ C}_z^* 
\times \displaystyle{ (\frac{ 1}w}\, \widetilde {D_b})  
& \ \ \mbox{if } \ [0,w] \in D\cap  {\bf  T}_b .
\end{array}\end{equation*}

We note the following:
 \begin{enumerate}
 \item [(1)] 
If $e\in D$ then $D[e]=\widetilde D \setminus \{zw=0\}= \widetilde {D^*}$; and 
$[z,w]\in D$ if and only if  $e \in D[z,w] $ (recall the definition of $\widetilde D$ (and hence $ \widetilde {D^*}$) in (\ref{tilded})); 
 \item [(2)] For each $[z,w]\in D$, $D[z,w]$ is an 
       open set with $C^\omega$ boundary $\partial D[z,w]$ but it is not relatively compact in $\mathbb{C}^* \times 
       \mathbb{C}^* $. We have  
\begin{enumerate}
 \item [(i)] $ D[z,w]= D[z,w] \times {\mathcal  I}$; 
 \item [(ii)] For $[z,w]\in D^*$ we define 
\begin{align*}
 D^*[z,w]&= \{(\alpha ,\beta )\in \mathbb{C}^* \times \mathbb{C}^*
: (\alpha ,\beta )[z,w]\in D^*\}.
\end{align*}
Then $D[z,w]=D^*[z,w]$. 
\end{enumerate}
 \item [(3)] 
\begin{enumerate}
 \item [(i)] 
For $[z,w]\in D^*$ we have 
\begin{align}
 \label{elm}
 D[z,w]= \widetilde {D^*} \times (\frac{ 
       1}{z}, \frac{ 1}{w})\ , 
\end{align}
and 
for $[z,w], [z',w'] \in D^*$ 
\begin{align} 
{  \label{parallel}}
D[z',w']&= (\frac{ z}{z'}, \frac{ w}{w'}) D[z,w].
\end{align}
In particular, the sets $D[z,w]$ for $[z,w]\in D^* $ are biholomophic in $\mathbb{C}^* 
 \times \mathbb{C}^* $.
 \item [(ii)] 
For any two points $[z,0], [z',0] \in D\cap {\bf  T}_a$
$$
D[z',0]= (\frac{ z}{z'}, 1) D[z,0].
$$
In particular, the sets $D[z,0]$ for $[z,0]\in  D\cap {\bf  T}_a $ are biholomophic in $\mathbb{C}^* \times \mathbb{C}^* $.
\end{enumerate}
\item [(3)] Fix  $[z_0,0]\in D \cap {\bf  T}_a$ and let $[z_n, w_n] \in  
D^* \ (n=1,2,\ldots )$ with $[z_n,w_n] \to [z_0,0]$ as $n\to 
       \infty$ in ${\mathcal  H}$.
For $0<r<R$, consider the product of annuli $$\mathcal A(r,R): \{r<|z|<R\}\times \{r<|w|<R\}\subset \mathbb{C}^* \times \mathbb{C}^*.$$
Then 
\begin{align}
 \label{approx}
\lim_{ n\to \infty} \partial D[z_n,w_n]\cap \mathcal A(r,R)= 
       \partial D[z_0, 0] \cap \mathcal A(r,R)  
\end{align}
in the Hausdorff metric as compact sets in $\mathbb{C}^* \times \mathbb{C}^* $. 

\end{enumerate}

We set 
\begin{align}
 \label{def-calD}
{\mathcal  D}:= \bigcup _{[z,w]\in D}\ ([z,w], D[z,w]).
\end{align}
This is a pseudoconvex domain in $ D \times (\mathbb{C}^* \times 
\mathbb{C}^*)$ which we consider as a function-theoretic ``parallel'' variation
$$
{\mathcal  D}: [z,w]\in D \to D[z,w] \subset \mathbb{C}^* \times 
\mathbb{C}^*. 
$$
Since $e\in D[z,w]$ for $[z,w]\in D$, we have the $c$-Green function 
$g([z,w], (\xi, \eta))$ with pole at $e$ and the $c$-Robin constant $\lambda[z,w] $ for 
$(D[z,w], e)$ with respect to the metric $ds^2$ on $\mathbb{C}^* \times 
\mathbb{C}^*$ and the function $c(z,w)>0$. 
We call $[z,w]\to \lambda [z,w]$ the {\it $c$-Robin function for $D$}. 

\vspace{2mm} 
The function $-\lambda [z,w]$ is a candidate to be a plurisubharmonic exhaustion function for 
 $D$. To be precise, we have the following fundamental result.
\begin{lemma} \label{lem:basic} 
 \quad \begin{enumerate}
 \item [1.]  $-\lambda [z,w]$ is a 
plurisubharmonic function on  $D$. 
 \item [2.] We have the following:
 \begin{enumerate}
 \item For any $[z_0,w_0]\in \partial D^*$, $\lim_{[z,w]\to [z_0,w_0]}\lambda [z,w]=-\infty$.
 \item If  $\emptyset \not = \partial D \cap {\bf  T}_a \not = {\bf  T}_a $ then for any 
$[z_0,0]\in \partial D \cap {\bf  T}_a$ we have $\lim_{[z,w]\to [z_0,0]}\lambda [z,w]=-\infty$ (and similarly if ${\bf  T}_a$ is replaced by ${\bf  T}_b$).
 \end{enumerate}
 \item [3.]  If $\partial D \not \supset {\bf  T}_a $  and 
$\partial D \not 
 \supset {\bf  T}_b$, then $-\lambda [z,w]$ is a plurisubharmonic exhaustion function for 
 $D$.
\end{enumerate}
\end{lemma} 
 
\noindent {\bf Proof.} Note that 3. follows from 1. and 2. We divide the proof of 1. into two steps.
\vspace{2mm} 

\noindent {\it  $1^{st}$ step.} \ \ {\it  $-\lambda [z,w]$ is plurisubharmonic on $D^*$}.
\vspace{2mm} 

\noindent Fix $[\zeta_0]=[z_0,w_0]\in D^*$. Let ${\bf  a}\in \mathbb{ C}^2 
\setminus \{0\}$ with $\|{\bf  a}\|=1$ and let $B=\{|t|<r\} \subset \mathbb{ 
C}_t$ be a small disk and let $(z(t),w(t))= \zeta_0 + {\bf  a}t$ be such that  the complex line 
 $l: t\in B
\to 
[\zeta (t)]=[z(t),w(t)] =[\zeta_0]+ {\bf  a}t $ passing through $[\zeta_0]$ is contained 
in $D^*$.
It suffices to prove that 
$-\lambda (t):= -\lambda [z(t), w(t)]$ is subharmonic on $B$, i.e., 
$$
\frac{\partial^2 \lambda(t)}{\partial t \partial 
\overline{ t}} \le 0\quad \mbox{ on $B$}.
$$
For brevity we write 
\begin{equation*}
\begin{array}{lll}
 & D(t): = D[\zeta (t)]  \subset \mathbb{C}^* \times \mathbb{C}^* 
&  \mbox{  for $t\in B$};\\[1mm]
&g(t,(z,w)) : = g([\zeta(t)], (z,w)) & \mbox{ for $(z,w) \in D[\zeta(t)]$}.
\end{array}
\end{equation*}
By (\ref{parallel}) we have 
\begin{align}
 \label{condition-1}
D(t)  =D[\zeta_0] \ (  \frac{ z_0}{z(t)}, 
\frac{w_0}{w(t)} ) \quad \mbox{ in $\mathbb{C}^* \times \mathbb{C}^* $.} 
\end{align}
 
We thus have the parallel variation of domains $D(t)$ in $\mathbb{C}^* \times 
\mathbb{C}^* $ with parameter $t\in B$:
$$
{\mathcal  D}|_{B}: t\in B \to D(t) \subset \mathbb{C}^* \times 
\mathbb{C}^* .
$$
We write 
 $${\mathcal  D}|_{B}:=\bigcup_{t\in B}\ (t, D(t)) ; \quad 
\partial {\mathcal  D}|_{B} = \bigcup_{t\in B}\ (t, \partial D(t))  \quad 
\mbox{ in } B \times (\mathbb{C}^* \times \mathbb{C}^*) ,
$$
where again we identify the variation with the total space ${\mathcal  
D}|_B$. By (\ref{def-calD}),   ${\mathcal  D}|_B$ is a pseudoconvex domain in 
$B \times (\mathbb{C}^* \times \mathbb{C}^*) $ such that  $\partial 
{\mathcal  D}|_B$ is $C^\omega $ smooth. Using the 
notation $\zeta=(z,w)\in \mathbb{C}^* \times \mathbb{C}^*    $ and  $g(t,\zeta)=g(t,(z,w))$,  
we have the following variation formula from Theorem 3.1 of \cite{kly}: 
\begin{align*}
(\star) \qquad  \frac{\partial^2 \lambda (t)}{\partial t \partial \overline{ t}}
= &-c_2 \int_{\partial D(t)}K_2(t,\zeta) \|\nabla _{\zeta}\ 
 g(t,\zeta)\|^2 dS_{\zeta}\\
&-4 c_2 \iint_{D(t)} 
\biggl( \bigl| \frac{\partial^2 g(t,\zeta)}{\partial \overline{  t} 
\partial z}\bigr|^2 + 
\bigl| \frac{\partial^2 g(t, \zeta)}{\partial \overline{  t} 
\partial w}\bigr|^2 \biggr) dV_{\zeta}\\
&-2c_2  
 \iint_{ D(t)}c(\zeta)\bigl|\frac{\partial g(t,\zeta)}{\partial 
t}\bigr|^2 dV_{\zeta}.
\end{align*}
Here $1/c_2$ is the surface area of the unit sphere in $\mathbb{ C}^2$; $dV_{\zeta}$ 
is the Euclidean volume element in $\mathbb{ C}^2$;  
\begin{align*}
 K_2(t,\zeta)&= {\mathcal  L}(t,\zeta)\bigl/ \|\nabla_{\zeta} \psi(t,\zeta)\|^{3}
\end{align*}
where 
${\mathcal L}(t,\zeta)$ is the ``diagonal'' Levi form defined by 
\begin{align*}
 {\mathcal  L}(t,\zeta) &= 
\frac{\partial^2 \psi}{\partial t \partial \overline{ t}} \, \|\nabla_{\zeta} \psi\|^{2} 
- 2 \Re\, \left\{
 \frac{\partial \psi}{\partial t }
\bigl(\frac{\partial \psi}{\partial \overline{ z} }
\frac{\partial^2 \psi}{\partial \overline{ t}  \partial z} 
+
\frac{\partial \psi}{\partial \overline{ w} }
\frac{\partial^2 \psi}{\partial \overline{ t}  \partial w}
\bigr) \right\}
+
\|\frac{\partial^2 \psi}{\partial t}\|^2\, \Delta _{\zeta} \psi;
\end{align*}
and $\psi(t,\zeta)$ is a defining function of ${\mathcal  
D}|_B$. The quantity $K_2(t,\zeta)$ is independent of the defining function $\psi(t,\zeta)$ (cf., Chapter 3 of \cite{kly}). Since ${\mathcal  D}|_B$ is pseudoconvex in $B \times (\mathbb{C}^* \times 
\mathbb{C}^*) $, following  Theorem 3.2 of \cite{kly} we have $K_2(t,\zeta) \ge 0$ on $\partial 
{\mathcal  D}|_B$ and hence $\frac{\partial^2 \lambda (t)}{\partial t \partial 
\overline{ t}}\le 0$ on $B$, proving the first step. 

Since $c(z,w)>0 $ in $\mathbb{C}^* \times \mathbb{C}^* $, the variation 
formula immediately implies the following rigidity result which will be useful later (cf., Lemma 4.1 of \cite{kly}).
\begin{remark} \label{basic-fact}
If \ $\frac{\partial^2 \lambda }{\partial t \partial \overline{ t}} 
 (0)=0$, then  $\frac{\partial g}{\partial t}(0, (z,w) ) \equiv 0$ on $D(0)$, i.e., 
$$\frac{\partial g ([\zeta_0]+{\bf  a}t, (z,w))}{\partial t}\bigl|_{t=0}\ \equiv
 0 \ \hbox{on} \ D[\zeta_0].$$ 
\end{remark}

\noindent {\it $2^{nd}$ step.} \ \ {\it Plurisubharmonic extension of 
$-\lambda [z,w]$ 
to $D$.}  \ 

\vspace{2mm} 
\noindent  We fix a point 
of $D \cap [(T_a \times \{0\}) \cup (\{0\} \times T_b)] )$, e.g., $[z_0, 
0]$ with $z_0\ne 0$. Let $[z_n, w_n]\in D^* \ (n=1,2,\ldots )$ with 
 $[z_n, w_n] \to [z_0,0]$ as $n\to \infty$. By 
 (\ref{approx}) 
\begin{align*}
& \lim_{ n\to \infty}(g([z_n,w_n], (\alpha , \beta) )- g([z_0,0], (\alpha 
,\beta)))=0\\
& \quad \mbox{  uniformly 
for $ (\alpha,\beta)$ in $K \Subset D[z_0,0] \subset 
 \mathbb{C}^* \times \mathbb{C}^* $.}
\end{align*}
It follows that $\lim_{ n\to \infty}\lambda [z_n,w_n]=\lambda [z_0,0]$, i.e.,
 $\lambda [z,w]$ is continuous and finite at $[z_0, 0]$. 
Hence $\lambda [z,w]$ is continuous and finite-valued on $D$. 
Since $D \cap {\bf  T}_a$ is a complex line, it follows from the first step 
that $-\lambda [z,w]$ extends to be plurisubharmonic from $D^* \cap {\bf  T}_a$ to $D \cap {\bf  T}_a$. Hence $-\lambda [z,w]$ extends to be plurisubharmonic on $D$. \hfill $\Box$ 

\vspace{3mm} 
We divide the proof of 2. in two steps; the first step is 2 (a).  
\vspace{2mm} 

\noindent {\it $1^{st}$ step.}  \ {\it Fix $[z',w']\in 
\partial D^*$. If $[z,w] \in D \to [z',w']$ in ${\mathcal  H}$,
then $\lambda [z,w] \to -\infty$.}

\vspace{2mm} 
\noindent Since $[z',w']\in \partial D^*$, we have $z'\ne 0$ and $w'
\ne 0$. 
If $[z,w]\in D^*$ tends to $[z', w']$ in ${\mathcal  H}$, 
then  $\partial D[z,w] 
\subset \mathbb{ C}^* \times  \mathbb{ C}^* $ tends to the single point 
$e$ in the sense that if we define $d [z,w]={\rm dist} (\partial D[z,w], 
e)>0$, where 
$$
{\rm dist} (\partial D[z,w], 
e):= {\rm Min} \ \{\sqrt{|\xi-1|^2+|\eta-1\}^2}:\ (\xi, \eta) \in 
\partial D[z,w]\},
$$
then $d[z, w] 
\to 0 $ as $[z,w]\to [z', w']$. Indeed, let $[z, w]\in D$ approach $[z', w']$ 
in ${\mathcal  H}$. By slightly deforming the fundamental 
domain ${\mathcal  F} \subset \mathbb{C}^* \times \mathbb{C}^*$ if necessary,  
we may assume $(z',w'), (z,w) \in {\mathcal F}$. Since 
$$
\partial D[z, w]= \{ (\frac{ \alpha }{z}, \frac{ \beta }{w}) 
\in \mathbb{ C}^* \times  \mathbb{ C}^* : \ [\alpha , \beta ]\in 
\partial D\}
$$
and $[z',w'] \in \partial D^*$, 
$$
d[z,w]= 
 {\rm dist}(\partial D[z,w], e)  
\le \sqrt{|{ z'}/{z}\, -1|^2+ |{ w'}/{w}\, -1|^2}$$
which clearly tends to $0$ as $[z,w] 
 \to [z',w']$. Since $\partial D[z,w]$ is a 
{\it smooth} real three-dimensional hypersurface, 
it follows by standard potential-theoretic arguments that 
$-\lambda [z,w] \to +\infty$. \hfill $\Box$ 


\vspace{2mm} 

It remains to prove 2 (b). Thus we assume $\emptyset \not = \partial D \cap {\bf  T}_a \not = {\bf  T}_a $.
\vspace{2mm} 

\noindent {\it  $2^{nd}$ step.}  \ {\it 
Fix $[z_0,0] \in \partial 
D \cap {\bf  T}_a$. If $[z,w]\in D \to [z_0,0]$ in ${\mathcal  H}$, then 
$ \lambda [z,w] \to -\infty$.}

\vspace{2mm}

For the proof of this step we require Lemma \ref{ps-lemma}. Fix $p_0=[z_0,0] \in \partial 
D \cap {\bf  T}_a$. We want to show
$$
\lim_{ [z,w] \to 
[z_0,0], \ [z,w]\in D } \lambda [z,w]=-\infty.
$$
We take a sequence $\{[z_n,w_n]\}_n\subset D$ which converges to $p_0$ in ${\mathcal  H}$. We show
\begin{align}
 \label{4thstep}
\lim_{ n \to \infty} \ \lambda [z_n,w_n]=-\infty.
\end{align}

From continuity of $\lambda [z,w]$ in $D$, it suffices to prove 
(\ref{4thstep}) for $[z_n,w_n] \in D^*$. Moreover, since $\partial D[z_n, w_n]$ is smooth, as in the end of the first step, we need only 
show
\begin{align}
 \label{distance-0}
\lim_{n \to \infty} {\rm dist} \ (\partial D[z_n,w_n], e )=0.
\end{align}

Before proving (\ref{distance-0}), we offer an example to explain the subtlety of the problem. We encourage the reader to draw a picture to illustrate the following situation. 
Let $D$ be a domain in ${\mathcal  H}$ with smooth boundary but which is not pseudoconvex. Precisely, we assume $D$ has the property that $ \partial D \cap {\bf  T}_a$ is a smooth curve in ${\bf  
 T}_a$ passing through a point  $[z_0,0] $ where $1<|z_0|<|a|$. 
We can find a bidisk $\delta := \delta _1\times \delta _2
 =\{|z-z_0|<r_1\} \times \{|w|< r_2\}\subset {\mathcal  F}$ with $r_1,r_2$ sufficiently small
so that $ D_1:=D \cap \delta $ is of the form 
$D_1= \cup _{z\in \delta _1}(z, D_1(z))$ where $D_1(z) \subset \delta 
_2$ and  $\partial D_1(z)$ is a 
non-empty smooth arc in $\delta _2$. 
{\it  We assume that, for each $z\in \delta _1$ 
\begin{align*}
 D_1(z) \supset D_1(z_0)  \supset  \delta _2\cap \{\Re\, w>0\} =: \delta_2^*.
\end{align*}} 
We can find a sequence $\{(z_n, u_n)\}_{n}$ in $D_1$ with $u_n=\Re\, w_n>0$ which converges to the point $(z_0, 0)\in \partial D$. Fix $r_1':0<r_1'< r_1/|z_0|$. 
By definition 
\begin{align*}
 D[z_n, u_n]&= (1/z_n, 1/u_n)\, \widetilde { D^*}  \supset (1/z_n, 
 1/u_n)\, \delta_1 \times  \delta_2^*
\end{align*}
and for sufficiently large $n$, say $n \ge n_0$,
$$E:= \{(Z,W)\in \mathbb{ C}^* \times \mathbb{ C}^*:|Z-1|<r_1', \ |W-1|< 
1/2\}\subset (1/z_n, 1/u_n)\, \delta_1 \times  \delta_2^*.$$
If we let $A$ denote the $c$-Robin constant 
for the domain $E$ in $\mathbb{ C}^* \times \mathbb{ C}^*$
and the point $e=(1,1)$, it follows that $\lambda [z_n, u_n]> A$ for $n   
\ge n_0$, so that ${  -} \lambda [z,w]$ is not an exhaustion function for $D$.

Returning to the proof of (\ref{distance-0}), we will use Lemma \ref{ps-lemma} and the {\it  pseudoconvexity} of the   domain $D$ in ${\mathcal  H}$. We may assume that $p_0=[z_0,0] \in \partial D$ lies in the fundamental domain ${\mathcal  F}$ and we take a sufficiently small bidisk $\Delta = \Delta _1 \times \Delta _2$ with  
center $(z_0,0)$ so that  $\Delta \subset {\mathcal  F}$. Let $\psi(z,w)$ be a defining 
function of $D$ in $\Delta $, i.e., $\psi(z,w)\in 
C^\omega(\Delta )$ with $D \cap \Delta =\{\psi(z,w)<0\}$ and $
\partial D \cap \Delta=\{\psi(z,w)=0\}$. Since $\partial D$ is smooth in 
${\mathcal  H}$, we have two cases: 
$$
 {\rm Case\  {  (c1)}}: \quad \frac{\partial \psi}{\partial z}\not \equiv 0 
 \mbox{ on} \  \Delta ; \qquad  {\rm Case\  {  (c2)}}: \quad \frac{ \partial \psi}{\partial w}\not \equiv 0 
      \mbox{  on} \  \Delta.
$$
Apriori, we also have two cases relating to the behavior of $\psi (z,0)$ on $\Delta_1$: 
$$
 {\rm Case\  { (d1)}}: \ \ 
\psi (z,0) \not \equiv 0 \ \ \mbox{ on $\Delta_1$}; \ \ \
 {\rm Case\  { (d2)}}: \ \ 
\psi (z,0) \equiv 0 \ \ \mbox{ on $\Delta_1$}.
$$
However, the hypothesis $\partial D \not \supset {\bf  T}_a$ 
in 2 (b) together with the 
real-analyticity of $\partial D$ imply that Case (d2) does not occur. Thus it suffices to prove (\ref{distance-0}) assuming that $\psi (z,0) \not \equiv 0$ on $\Delta_1$.

\vspace{2mm} 
{{\it    Proof of (\ref{distance-0}) in Case (c1):} 

\vspace{2mm} 
In this case, by taking a suitably smaller bidisk $\Delta $ if necessary, $l(0):=\{\psi(z,0)=0\}$ is a $C^\omega-$smooth arc in $\Delta _1$ 
passing through $z=z_0$ and $l(0) \times \{0\} \subset \partial D\cap \Delta $. For $w\in \Delta _2$, 
$$
l(w):= \{ z\in \Delta _1: (z,w) \in \partial D \cap \Delta \}.
$$
is a simple $C^\omega-$smooth arc in $\Delta _1$.

Fix $\varepsilon >0$. Since $z_0\ne 0$, we can find a disk $\delta _1 \subset \Delta _1$ 
with center $z_0$ 
such that  
$$
|\frac{ z'}{z''}-1|<\varepsilon \quad \mbox{ for all $z',z''\in \delta _1$}.
$$
Now we take $\delta _2: |w|<r <\varepsilon $ in $\Delta _2$ so that  
each arc $l(w)$ passes through a certain point $\zeta (w)$ in $\delta_1$. For sufficiently large $n_0$, if $n \ge n_0$ we have $(z_n, w_n) \in \delta _1 \times \delta _2$. Since 
 $w_n\in \delta _2$, we have $\zeta (w_n)\in l(w_n) \cap \delta_1$ so that $(\zeta (w_n),w_n) 
\in \partial D$ in ${\mathcal  H}$. Hence, $(\frac{ \zeta (w_n)}{z_n}, \frac{ w_n}{w_n})=
(\frac{ \zeta (w_n)}{z_n}, 1) \in \partial D[z_n,w_n]$
in  $\mathbb{ C}^* \times  \mathbb{ C}^* $. Thus
$$
{\rm dist}\ (\partial D[z_n,w_n],e\ ) \le 
 {\rm dist}\ ((\frac{ \zeta (w_n)}{z_n}, 1),e\ )= |\frac{ \zeta (w_n)}{z_n}-1| < 
 \varepsilon \quad   \ \hbox{for} \ n \ge n_0.
$$

{\it  Proof of (\ref{distance-0}) in Case ({c2}):} 

\vspace{2mm} 
In this case, by taking a suitably smaller bidisk $\Delta $ if necessary, 
$S(z_0):=\{\psi(z_0,w)=0\}$ is a $C^\omega-$smooth arc in $\Delta _1$ 
passing through $w=0$ and $ \{z_0\} \times S(z_0) \subset \partial D\cap \Delta $. For $z\in \Delta _1$, 
$$
S(z):= \{ w\in \Delta _2: (z,w) \in \partial D \cap \Delta \},
$$
is a simple $C^\omega-$smooth arc in $\Delta _2$.

Fix $\delta _1:=\{|z-z_0|<r_1\} \Subset D(z_0)$. 
Case (d1) corresponds to the condition ${(iii)}$ in Lemma 
\ref{ps-lemma}, thus this lemma implies that there exists a disk $ \delta _2:= 
\{|w|<r_2\} $  such that  
\begin{align}
  \label{letitbe-3}
 \bigcup_{z\in \delta_1} S(z)\ {  \supset  }\
D(z_0) \cap\delta_2.
\end{align}

Fix $\varepsilon >0$. Taking $r_1$ sufficiently small, we can insure that $$
|{ z'}/{z''}\, -1|< \varepsilon  \quad \hbox{ for all} \ \ z',\ z''\in \delta _1.
$$
Take a disk $\delta _2 \subset \Delta _2$ satisfying (\ref{letitbe-3}).  For sufficiently large $n_0$, if $n \ge n_0$ we have $(z_n, w_n) \in \delta _1 \times \delta _2$. 
We divide the points $w_n\in \delta _2$ into two types:

$$
{\rm Case \ (i):}\  \ w_n \in \delta _2 \cap D(z_0); \ \ \ {\rm Case \
(ii)}:  \ w_n \in \delta 
_2 \setminus D(z_0).
$$

In Case (i), using (\ref{letitbe-3}) we can find $z^*\in \delta _1$ with $w_n \in S(z^*)$  
so that  $(z^*, w_n)\in \partial D$ in ${\mathcal  H}$ 
(see { $w_n, \ z^*, \partial D(z^*)$} in the figure below). Thus, 
 $({ z^*}/{z_n}, { w_n}/{w_n})=({ z^*}/{z_n}, 1)$ in 
$ \partial D[z_n,w_n]$ in 
 $\mathbb{ C}^* \times  \mathbb{ C}^* $ and hence 
$$
{\rm dist}\ (\partial D[z_n,w_n],e\ ) \le 
 {\rm dist}\ ({ z^*}/{z_n}, 1),e\ )= |{ z^*}/{z_n}\, -1| < 
 \varepsilon \quad  \hbox{ for all} \  \ n \ge n_0.
$$

In Case (ii), let $\ell=[z_n, z_0]$ be a segment in $\delta _1$. We can find $z^*\in \ell $ with $w_n \in \partial D(z^*)$. Indeed, as $z$ goes from $z_n$ to $z_0$ along $\ell$, 
the arcs $\partial D(z)\cap \delta _2$ transform from $\partial D(z_n)\cap \delta _2$ to $\partial D(z_0)\cap 
\delta _2$ in a continuous fashion. Since $[z_n, w_n] \in D^*$, we can find $z^*\in \ell $ with $w_n \in \partial D(z^*)$. 

Thus $(z^*, w_n)\in \partial D^*$, so that
$(z^*/z_n, 1) \in \partial D^*[z_n,w_n]$, and hence
$$
{\rm dist}\ (\partial D[z_n,w_n],e\ ) \le 
 {\rm dist}\ ({ z^*}/{z_n}, 1),e\ )= |{ z^*}/{z_n}\, -1| < 
 \varepsilon \ \  \hbox{ for all}\ n \ge n_0,
$$
which is (\ref{distance-0}). This completes the  proof of 2 (b) in Lemma \ref{lem:basic}.\hfill $\Box$ 
\vspace{2mm}

 We next relate the possible absence of strict plurisubharmonicity of the function $-\lambda [z,w]$ on a pseudoconvex domain $D$ in ${\mathcal  H}$ {\it at a point in $D^*$} with existence of holomorphic vector fields on ${\mathcal  H}$ with certain properties. 
This is in the spirit of, but does not follow from, Lemma 5.2 of \cite{kly}. Recall that in the case $\rho$ is rational and $\tau$ is rational, we defined $ \sigma_c:= \{w=c z^\rho\}/ \sim $ to be the integral curve $[z_0,w_0] 
       \exp tX_u$ with $c= w_0/z_0^\rho\ne 0, \infty$ of $X_{u}: = (\log|a|)\, z\, \frac{\partial }{\partial z} + (\log|b|) \,
w\,\frac{\partial }{\partial w}$.

\begin{lemma}
 \label{essen-lemm} 
Let $D$ be a pseudoconvex domain with $C^\omega-$smooth boundary in 
 ${\mathcal  H}$ and let $\lambda [z,w]$ be the $c$-Robin function on $D$.
 Assume that there exists a point $p_0= [z_0, w_0]$ in $D^*$ at which $-\lambda 
       [z,w]$ is not strictly plurisubharmonic.  
\begin{enumerate}
 \item [(1)]  There exists a 
holomorphic vector field $X=\alpha z \frac{\partial }{\partial z} dz 
+\beta w\frac{\partial }{\partial w }dw \ne 0$ on ${\mathcal  H}$ such that 
if $[z,w]\in D^* \ (resp. \ \partial D^*)$, then the integral curve 
       $I[z,w]:= [z,w]\exp tX$  in 
       ${\mathcal  H}$ is contained  in $D^* \ (resp.\ \partial D^*)$. We say $X$ is a {\sl tangential vector field on} $\partial D^*$.
 \item [(2)] The form of the vector field $X$ in (1) and the domain $D$ are determined as follows:
\begin{enumerate}
 \item [(i)]
\ if $ \partial D\not \supset {\bf  T}_a$ and 
       $\partial D\not \supset {\bf  T}_b$, 
then  $X =cX_u$ for some $c\ne0$ with $X_u$ in (\ref{xu}). 
If $\rho$ is irrational or $\rho$ is rational and $\tau$ is irrational, $D$ is of 
type {\rm (a-1)} or of type {\rm (b1)} in Theorem 
\ref{thm-1}. If $\rho$ is rational and $\tau$ is rational, $D$ is of 
       type {\rm (b2)} :\   
$D=\cup _{c\in \delta}\, \sigma_c$ where 
       $\delta $ is a relatively compact domain in $\mathbb{ C}^*$ 
with smooth boundary.
 In all cases, we have $\partial D \cap ({\bf  T}_a \cup {\bf  T}_b) 
       =\emptyset$. 
 \item [(ii)]  \ if $\partial D \supset {\bf  T}_a $ and $\partial D 
       \not \supset  {\bf  T}_b$, 
then we have two cases:
\begin{enumerate}
\item[(ii-a)]  $X= cX_u $ for some $c\ne 0$ and $D$ is of type 
	     ${\rm (b2)}$: $D=\cup _{c\in \delta}\, \sigma_c$ where 
       $\delta $ is a domain in $\mathbb{ P}^1=\CC\cup \{\infty\}$ with smooth boundary $\partial\delta 
       $ which contains $0$ but not $\infty$.  
 \item [(ii-b)] $X=c z \frac{\partial }{\partial z}$ 
 for some $c \ne 0$. Then $D$ is a domain of ``Nemirovskii 
	     type'':
  $b >1$  and $D= \mathbb{ C}_z \times \{Au+ B v<0\}/\sim $, where              
	     $A,B\in \mathbb{ R}$ with $(A,B)\ne (0,0)$ (here $w=u+iv$).
\end{enumerate}
\item [(ii')]  if  $\partial D \supset {\bf  T}_b $ and $\partial D 
       \not \supset  {\bf  T}_a$, we have the result analogous to (ii).
 \item [(iii)] \ If $\partial D \supset {\bf  T}_a \cup {\bf T}_b$, then 
$X= cX_u $ for some $c$ and $D$ is of type {\rm (b2)}:\  
$D=\cup _{c\in \delta}\, \sigma_c$ where 
       $\delta $ is a domain in $\mathbb{ P}^1$ with smooth boundary $\partial\delta 
       $ with $0, \infty \in \partial \delta$.
\end{enumerate}
\end{enumerate}
\end{lemma}
\begin{remark}\label{nemrem} With respect to the Nemirovskii-type domain in (ii-b), we recall Nemiroviskii's theorem in \cite{nem}. Let $a>1$ and let ${\mathcal H}={\mathcal  H}_{a,a} $. Then the domain  
$D=\mathbb{ C}_z \times \{\Re\, w>0\} / \sim \  \subset {\mathcal  H}$ is Stein and 
$\partial D$ is Levi-flat. At the end of section 4 we will discuss an explicit example of such a domain which will illustrate some of the ideas used in the proof of Theorem \ref{thm-1}.
\end{remark}
\noindent {\bf  Proof.} \ Since $\lambda [z,w]$ is plurisubharmonic on $D$ 
and is not strictly plurisubharmonic at $p_0=[z_0,w_0]\in D^*$, we can find a 
holomorphic vector field 
$
X= \alpha z \frac{\partial }{\partial z}dz + 
\beta w \frac{\partial }{\partial w}dw \ne 0$ on 
${\mathcal  H}$ such that  
\begin{align}\label{condition-n}
 \frac{\partial^2 \lambda [\ p_0\exp t X\ ]}{\partial t \partial 
\overline{ t}} \ \biggl|_{t=0} =0.
\end{align}
We shall show that this $X$ is a tangential vector field on $\partial D^*$. Since $p_0\in D^* $, we can take a small disk $B=\{|t|<r\}$ 
with $p_0 \exp tX 
\subset D^*$ for $t\in B$. 
We set $D(t)=D[p_0 \exp tX] \subset \mathbb{C}^* \times 
\mathbb{C}^* $ so that $D(0)=D[p_0]$. We let $g(t, (z,w))$ (resp. $\lambda (t)$) denote the 
$c$-Green function $g([p_0 \exp tX], (z,w) )$ (resp. the $c$-Robin 
constant $\lambda[p_0 \exp tX]$) for $(D(t), e)$ 
and $t\in B$. We set
${\mathcal  D}|_B= \cup_{t\in B}(t, D(t)) \subset B 
\times (\mathbb{C}^* \times \mathbb{C}^* )$ which we consider as the variation
$$
{\mathcal  D}|_B: \ t\in B \to D(t)=D[p_0 \exp tX] \subset \mathbb{C}^* \times \mathbb{C}^*.
$$
By  
(\ref{parallel}) we have 
\begin{align*}
D(t)&= D[p_0 \exp tX]= D[[z_0,w_0]\exp tX] \\
&= D[z_0, w_0]  \exp (-tX) = D[z_0,w_0] \bigl( e^{-\alpha t}, 
  e^{-\beta  t}\bigr)  \quad \mbox{ in $\mathbb{C}^* \times \mathbb{C}^* $.}
\end{align*}
Using the same reasoning as in the first step of the proof of Lemma \ref{lem:basic} 
together with Remark \ref{basic-fact} we see from (\ref{condition-n}) and the real analyticity of $\partial {\mathcal  D}|_B= 
\cup_{t\in B}(t, \partial D(t))$ in $B \times (\mathbb{C}^* \times \mathbb{C}^*)$ that 
\begin{align}
\label{norm}
\frac{\partial g(t,(z,w)}{\partial t}\biggl|_{t=0} { \equiv}\ 0 \quad 
\mbox{ on $D[z_0,w_0] \cup \partial D[z_0,w_0]$.}
\end{align}
For a fixed $t\in B$ we consider the automorphism
$$
(Z,W) \to (z,w)=F(t, (Z,W))
$$
 of $\mathbb{C}^* 
\times \mathbb{C}^* $ where 
\begin{align*}
 F(t, (Z,W)):=
 (Z,W)\bigl( \frac{ 
1}{z_0}, \frac{ 1}{w_0}\bigr) \exp (-tX)=
\bigl(\frac{Z e ^{-\alpha t}}{z_0}, \frac{W  e^{-\beta 
 t}}{w_0} \bigr). 
\end{align*}
Then
$$
(z,w)\to (Z,W)=F^{-1}(t, {  (z,w)})= \biggl (zz_0e^{\alpha t}, ww_0e^{\beta 
t} \biggr).
$$
By (\ref{elm}) we have 
 \begin{align*}
  D(t)&= \widetilde {D^*} \,
 \bigl( \frac{ 
1}{z_0}, \frac{ 1}{w_0}\bigr) \exp (-tX) = \widetilde {D^*}\, \bigl(\frac{ e ^{-\alpha t}}{z_0}, \frac{  e^{-\beta 
 t}}{w_0} \bigr)  \quad \mbox{ in $\mathbb{C}^* \times \mathbb{C}^* $,}
\end{align*}
so that 
 $D(t)= F(t, \widetilde D^*)$. We   note that 
$\widetilde {D^*} \subset \mathbb{C}^* \times \mathbb{C}^* $ is independent of $t\in B$.
 We set 
$$
G(t,(Z,W)):= g(t, (z,w)) \quad \mbox{ where $(z,w)=F(t,(Z,W)), \quad  
(Z,W)\in \widetilde { D^*}$}.
$$
Since
$$
g(t,(z,w))=G(t, F^{-1}(t,(z,w))= G(t, (zz_0 e^{\alpha t}, ww_0 e^{\beta t})),
$$
we have 
\begin{align*}
&\frac{\partial g} {\partial t}(t,(z,w))\\
&= 
\frac{\partial G }{\partial t}(t,(Z,W))
+\frac{\partial G}{\partial Z }(t,(Z,W) ) \alpha zz_0  e^{\alpha t} 
+  \frac{\partial G}{\partial W }(t,(Z,W) ) \beta ww_0  e^{\beta t}\\
&=
\frac{\partial G }{\partial t}(t,(Z,W))
+ \alpha Z  \frac{\partial G}{\partial Z }(t,(Z,W) )
+  \beta W \frac{\partial G}{\partial W }(t,(Z,W) )
\end{align*}
where $(Z,W)= F^{-1}(t,(z,w))$.  Since, for each $t\in B$, 
\begin{align}
 \label{G-function}
G(t, (Z,W)) \equiv 0 \mbox{ \quad on $\partial \widetilde {D^*}$},
\end{align}
we have 
$$
\frac{\partial G}{\partial t}(t, (Z,W)) =0 \quad \mbox{ on $\partial 
\widetilde {D^*}$}.
$$
It follows from (\ref{norm}) that  
\begin{align*}
\alpha Z  \frac{\partial G}{\partial Z }( 0,(Z,W) ) 
+ \beta W \frac{\partial G}{\partial W }(0,(Z,W) )=0 
 \quad \mbox{ on $\partial \widetilde {D^*}$.}
\end{align*}
Together with (\ref{G-function}), this says that the holomorphic vector field 
$$
X= \alpha Z\frac{\partial }{\partial Z}+ 
\beta W \frac{\partial }{\partial W},
$$
considered as a vector field on on $\mathbb{C}^* \times \mathbb{C}^* $, satisfies the property that for any $(z,w)\in \partial \widetilde {D^*}$, the integral curve 
$(z,w) \exp tX \subset \partial \widetilde{D^*}$ for all $ t\in \mathbb{ C}$. It follows that for 
 any $(z,w) \in \widetilde {D^*}$, the integral curve 
$(z,w) \exp tX$ is contained in $\widetilde {D^*}$: 
$$
\widetilde {D^*} \exp tX = \widetilde {D^*}, \ \hbox{for all} \ t\in \mathbb{ C}.
$$
Hence $X$ is a tangential vector field on 
$\partial \widetilde {D^*}$.

This implies
\begin{align}
 \label{rc}
D[ [z,w]\exp tX] = D[z,w]\subset \mathbb{C}^* \times 
\mathbb{C}^*, \ \hbox{for all} \  t\in \mathbb{ C} \end{align}
if $[z,w]\in D^*$ since 
\begin{align*}
 D[ [z,w]\exp tX] &= \widetilde {D^*}\ (\frac{ 1}z, \frac{ 1}{w})\ \exp 
 (-tX) = \widetilde {D^*}\ (\frac{ 1}z, \frac{ 1}{w}) = D[z,w].
\end{align*}
But for {  $[z,w]\in D^* \  (\hbox{resp.}\ \partial D^*)$} it is 
clear that 
\begin{align*}
 &[z,w] \exp tX \subset D^* \ (\hbox{resp.}\ \partial D^*) \mbox{  in ${\mathcal 
 H}$} \\
& \qquad \qquad  \mbox{ if and only if} \\
&(z,w) \exp tX \subset \widetilde {D^*} \ 
(\hbox{resp.} \ \partial \widetilde {D^*}) \ \mbox{ in $\mathbb{C}^* \times 
\mathbb{C}^* $},
\end{align*}
which proves that $X$, as a holomorphic vector field on ${\mathcal H}$, is a tangential vector field on $\partial  {D^*}$, verifying (1) of Lemma \ref{essen-lemm}. 
\vspace{3mm} 

To prove assertion (2) we first observe by (\ref{rc}) 
$$
\lambda [z,w]= \lambda [[z,w]\exp tX], \ \hbox{for all} \ t\in \mathbb{ C}
$$
for any $[z,w] \in D^*$. In case (2)(i) in Lemma \ref{essen-lemm}, from 3 in Lemma \ref{lem:basic},   
the Robin function $-\lambda[z,w] $  is an exhaustion 
function on $D$, and it follows that  
\begin{equation}\label{firstfact}
\{[z,w] \exp tX: t\in \mathbb{ C}\} \Subset D \ \hbox{for} \ [z,w]\in D^*.
\end{equation}

We now prove (2) (i). First we show that $X=cX_u$ for some $c\ne0$. If not, i.e., if  
$X\not \in \{cX_u: c \in \mathbb{ C}^*\}$, we take $[z,w] \in \partial D^*$ and 
let $\sigma=[z,w]\exp t X $ be the integral 
curve of $X$ passing through $[z,w]$. From Lemma \ref{fund-1} part 2 (2), 
the closure $\Sigma $ of $\sigma$ in ${\mathcal   H}$ contains 
 ${\bf  T}_a $ or 
${\bf  T}_b$, which contradicts the hypothesis $ \partial D\not \supset {\bf  T}_a$ and 
       $\partial D\not \supset {\bf  T}_b$ of (2) (i) in Lemma \ref{essen-lemm}.
Thus $X=cX_u$ for some $c\ne 0$. 

By (\ref{firstfact}), for $[z,w]\in 
D^*$ the closure of the integral curve $I[z,w]:=[z,w]\exp tX_u$ 
is compactly contained in $D$ and hence lies in $D^*$. 
{ It follows from ($\alpha$) and ($\beta$) in section 3}   
that we have
\begin{align*}
 (\alpha^* )& \ \  D^* =\bigcup _{c\in I}\ \{|w|=c |z|^\rho\}{ /\sim},  \mbox{ where 
$I$ is an open interval in $(0, \infty)$; or}\\
(\beta^*)& \ \  D^* =\bigcup _{c\in \delta}\ \{w=c z^\rho\}{ /\sim}, \ \mbox{where  
      $\delta$ is a domain in $\mathbb{ C}^*$.}
\end{align*}

We next show that if $D\cap {\bf T}_a\not =\emptyset$ then $D\supset {\bf T}_a$, contradicting the hypothesis in (2) (i). Thus let $[z_0,0]\in D\cap  {\bf T}_a$. Let $U,V$ be sufficiently small disks such that
$$(z_0,0)\in U\times V =:U \times \{|w|<r\}\Subset D \cap E_2$$
where recall $E_2= \{1 {  <} |z|\le |a|\} \times \{ |w|\le |b|\}\subset {\mathcal F}$. 
We show that there exists $r'$ with $0<r'<r$ such that 
\begin{equation}\label{alphabeta}G:=\{(z,w)\in E_2:1<|z|<a, \  0<|w|<r'\}\subset D.\end{equation}

We set }
$$  \Delta :=\{c=w/z^\rho \in \mathbb{ C}^*: (z,w)\in U \times \{0<|w|<r\},
$$
 so that 
$\Delta $ contains a punctured disk 
$\delta':= \{0<|c|<r'\} $ in $\mathbb{ C}$. Here we can take, e.g., $r'=r/|z_0|$ and we show this $r'$ works to achieve (\ref{alphabeta}).  
Using (1) in Lemma \ref{essen-lemm} we have $[z,w] \exp t X_u 
\subset D$ for $[z,w] \in U \times \{0<|w|<r\}$. 
To verify (\ref{alphabeta}),  
given  $ [z_0,w_0]\in G$, we have $c_0:= 
w_0/z_0^\rho \in \delta ' \subset \Delta$. Thus we can find $(z_0',w_0') \in U \times \{0<|w|<r\}$ 
with $c_0=w'_0/(z'_0)^\rho$. It follows that
$$[z_0, w_0] \exp t X_u = \{w= c_0z^\rho\}/\sim \ = [z_0', w_0']\exp tX_u 
\subset D.$$
In particular,  $[z_0, w_0]\in D$; hence (\ref{alphabeta}) is proved.

Suppose $D\not \supset {\bf T}_a$. We use the pseudoconvexity of $D$ to derive a contradiction. Observe that $D(0):=D\cap {\bf T}_a$ is a domain in ${\bf T}_a$ whose boundary $l$ consists of  smooth real one-dimensional curves. Fix $z'\in D(0)$ near $l$. Let $D(w)$ denote the slice of $D$
corresponding to $w$ for $0 < |w| < r'$. We consider the Hartogs radius $r(w)$ for $D(w)$ centered at $z'$. By (\ref{alphabeta}), $r(0) < r(w)$ for $0 < |w| < r'$. Since $D \cap E_2$ is pseudoconvex in $E_2$, this contradicts the superharmonicity of $r(w)$. A completely similar argument shows that if $D\cap {\bf T}_b\not =\emptyset$ then $D\supset {\bf T}_b$. Thus either $D=D^*$ as in ($\alpha^*$) or ($\beta^*$) or $D\setminus D^*$ consists of ${\bf T}_a, \ {\bf T}_b$, or ${\bf T}_a\cup {\bf T}_b$ with $D^*$ as in ($\alpha^*$) or ($\beta^*$). We verify that $D \setminus D^*={\bf  T}_a$ cannot happen; entirely similar proofs show that $D \setminus D^*={\bf  T}_b$ and $D \setminus D^*={\bf  T}_a\cup {\bf T}_b$ cannot occur. Indeed, if 
$D \setminus D^*={\bf  T}_a$, then $\partial D= {\bf  T}_a$, which is 
a complex line. However, $\partial D$ is assumed to be smooth; hence it must be a  
real three-dimensional surface. This completes the proof of (2) (i).


To prove (2) (ii), we note that under the 
 condition $\partial D \supset {\bf  T}_a $ and $\partial D 
       \not \supset  {\bf  T}_b$, from Lemma \ref{fund-1} we have either $X= cX_u$ with $c\ne 0$ or $X= \alpha z \frac{\partial 
 }{\partial z} $ with $\alpha \ne 0$. Using the same reasoning as in the proof of 2 (i) we conclude that $D$ cannot be of the form 
in Case {\bf a} nor of the form in Case ({\rm b1}) in Theorem \ref{thm-1}. 
 
 If $X= cX_u$ with $c\ne 0$, then $D^*$ is of the form $(\beta^*)$. Since $\partial D \supset {\bf  T}_a $ and $\partial D 
       \not \supset  {\bf  T}_b$ we arrive at the conclusion in (2) (ii-a). On the other hand, if $X= \alpha z \frac{\partial 
 }{\partial z} $ with $\alpha \ne 0$, we first observe from the facts that $\partial D \supset {\bf  T}_a$ and 
 $\partial D$ is $C^\omega-$smooth, for any $z_0\in 
 \mathbb{ C}^* $ the slice of $\partial D$ over $z=z_0$ {  contains} 
 a $C^\omega$ 
curve $C(z_0) \subset \mathbb{ C}_w$  passing through the origin $w=0$.
We can find a sufficiently small disk $V:= \{|w|< r_0\}$ so that $C(z_0)$ divides 
$V$ into two parts $V'$ and $V''$ with $\{z_0\} \times V' \subset D$ and $\{z_0 \} 
\times V''\subset \overline{ D}^c$. We set $\widetilde C(z_0):=C(z_0)\cap V$. 
By (1) in Lemma \ref{essen-lemm} 
we conclude that $\mathbb{ C}^* \times V' \subset D$ and $\mathbb{  C}^* \times 
V'' \subset \overline{ D}^c$. Thus $\mathbb{ C}^* \times \widetilde C(z_0)
 \subset 
\partial D$, which implies $\partial D \cap (\mathbb{ C}^* 
\times V) = 
\mathbb{ C}^* \times \widetilde C(z_0) $ and $D \cap (\mathbb{ C}^* \times V)= 
\mathbb{ C}^* \times V'$. 

We use this geometric set-up to show that $b$ must be a positive real number 
(hence $b>1$). To see this, fix a point $w_0\in \widetilde C(z_0)$ (resp.\,$V'$) with $w_0 \ne 0$. Since  $(z_0, 
w_0)\in \partial D$ (resp. $V'$), we have $\mathbb{ C}^* \times \{w_0\} \subset 
\partial D$ (resp. $D$). In particular, $(a^nz_0, w_0) \in 
\partial D$ (resp. $D$) for any $n \in \mathbb{ Z}$. Hence $(z_0, 
w_0/b^n) \in \partial D$ (resp. $D$) for any $n \in \mathbb{ Z}$.  
Since $|b|>1$ we can take $N$ sufficiently large so that $w_0/b^N\in V$. It follows that 
$ w_0/ b^n \in \widetilde C(z_0)$ (resp. $V'$) 
for any $n \ge N$. 

  We first show that $b$ is real. If not, let $b=|b|e^{i\phi}$ where $|b|>1$ and 
$ 0<|\phi|<\pi$. We set $w_0=|w_0|e^{i \varphi 
_0}$. Let  ${\mathbf  
n}_0=e^{i\theta_0}$ be a unit normal vector to $\widetilde C(z_0)$ at 
$w=0$ pointing in to $V''$.  Since $\widetilde C(z_0)$ is smooth, we can find $r_1$ sufficiently small 
with $0<r_1< r_0$ so that the sector ${\mathbf  e}:= \{r e^{i\theta}: 0<r<r_1, \ 
|\theta-\theta_0|< 2\pi/3\}$ is contained in $V''$. For any $N'\in \mathbb{ Z}$, 
it is clear that there exists $n'> N'$ satisfying 
\begin{align}
 \label{positive}
 \ |(\varphi _0- n'\phi) - \theta_0| < 2\pi/3 \ \hbox{modulo} \ 2\pi.
\end{align}
We take $N'>N$ so that $|w_0|/|b|^{N'} < r_1$, and then we 
choose $n'>N'$ with property (\ref{positive}). Then $w_0/b^{n'} \in {\bf  e} 
\subset V''$, which contradicts the fact that $w_0/ b^{n'} \in \widetilde C(z_0)$. Thus $b$ is real.

We next show $b$ is positive. If not, we have $b<-1$. We take $w_1\in 
V'\setminus \{0\}$ close to $0$. Then $(z, w_1) \in  D $ for all $z\in 
 \mathbb{ C}^*$. In 
particular, $(a^n z_0, w_1) \in D$ for any 
 $n \in \mathbb{ Z}$; hence $(z_0, w_1/ b^n) \in (\{z_0\} \times  
 V) \cap D $ for $n$ sufficiently large. In other words, for $n>N$ we have $w_1/b^n \in V'$. Since $b<-1$ it follows that $\{w_1/b^n: n \ge N\}$ lies on a line $L$ passing 
 through $w=0$. Moreover, if we take a sufficiently small disk $
 V_0:=\{|w|<r_0\} \subset V$, then $L \cap V_0$ intersects 
 the smooth curve $\widetilde C(z_0)$ transversally. At the point $
 w=0$, $L\cap V_0$ divides into 
 two segments $L'$ and $L''$ with $L'= (L\cap V_0) \cap V' $ and $L''=(L\cap V_0)\cap 
 V''$. Since $b<-1$, for $n$ sufficiently large, if $w_1/b^n \in L'$ then 
 $w_1/b^{n+1}\in L''$. This contradicts the fact that  
 $w_1/b^m\in V'$ for all $m$ sufficiently large. Thus $b>1$.

Consequently,  
$$
w \in \tilde C(z_0) \ ({\rm resp.} V') \longrightarrow { w}/{b^n} \in 
\tilde C(z_0) \ ({\rm resp.} V')  
\quad \mbox{ for $n=1,2,\ldots$. }
$$
It follows from the smoothness of $\tilde C(z_0)$ and the fact that $b>1$ that $\tilde C(z_0)$ is a line $Au+Bv=0$ passing through $w=0$, proving (2) (ii-b).

To verify (2) (iii), we show 
\begin{equation}\label{xin} X \in \{cX_u: c\in \mathbb{ C}\}\cup 
\{\alpha z\frac{\partial }{\partial z }: \alpha \in \mathbb{ C}\}\cup \{\beta w \frac{\partial }{\partial 
w}: \beta \in \mathbb{ C}\}.\end{equation}
Once (\ref{xin}) is verified, we obtain 2 (iii) by repeating 
the arguments in 2 (i) and 2 (ii). Suppose $X=\alpha z \frac{\partial }{\partial z}+ 
\beta w\frac{\partial }{\partial w} \not \in \{cX_u: c\in \mathbb{ C}\} 
 $ where $\alpha \ne 0, \beta \ne 0$. We set $\beta /\alpha =A+ iB$ where 
 $A, B$ are real numbers. To get a contradiction, we work in the case where $A$ is irrational; the case when $A$ is rational is similar. Fix $z_0 \in \{1< |z|< |a|\}$. 
 Since $\partial D \supset {\bf T}_a\cup {\bf  T}_b$ and $\partial D$ is 
  smooth, we can find a smooth curve $\ell$ in $\{|w|<|b|\}$ containing $w=0$ such 
  that $\{z_0\} \times \ell \subset \partial D$. We fix a disk 
  $V: =\{|w|<r\}$ with $r$ sufficiently small so that $\ell$ divides $V$ into two parts $V'$ and 
  $V''$ where $\{z_0\} \times  V'\subset D$ and $\{z_0\} \times V'' \subset 
  \overline{ D}^c$. Let $w_0 \in V'$ and for $c= w_0/z_0^\rho$, we consider the integral curve 
  $\sigma _c:= \{w= 
  c z^{A+iB}\} {  /\sim}$ of $X$ passing through $(z_0, w_0)$ 
  in ${\mathcal  H}$. Using (1) in 
  Lemma \ref{essen-lemm} we see that $\sigma _c \subset D$. On the other 
  hand, by Remark \ref{al-be-nonzero} there is a point $(z_0, 
  w(z_0))\in \sigma _c$ with $w(z_0) \in V''$, which is a 
 contradiction. This proves (\ref{xin}) and hence 2 (iii).  \hfill $\Box$

\vspace{3mm} 

Given a pseudoconvex domain $D$ in ${\mathcal  H}$ with $C^\omega-$smooth boundary, under the various cases of (2) of Lemma \ref{essen-lemm}, depending on the relationship between the tori ${\bf  T}_a, \ {\bf  T}_b$ and $\partial D$, we want to show that either $D$ is Stein or $D$ is the appropriate type of non-Stein domain in Theorem \ref{thm-1}. This will be done in a series of lemmas. Before proceeding, we recall an important ``rigidity'' result from \cite{kly}. 

We let ${\mathcal  D}: t\in B \to D(t) \subset M$ be 
a smooth variation of domains $D(t)\subset  M$ over $B\subset \CC$ where $M$ is a complex Lie group of dimension $n\ge 1$. Here $D(t)$ need not be relatively compact in $M$ but $\partial D(t)$ is assumed to be $C^\infty-$smooth. Assume each domain $D(t)$ contains the identity element $e$. Let  $g(t,z)$ and $\lambda (t)$ be the $c$-Green function and the $c$-Robin constant for $(D(t), e)$ associated to a K\"ahler metric and a positive, smooth function $c$ on $M$. We have the following from \cite{kly}:
\vspace{2mm} 
 
$(\star 1)$ \ \  {\it  Assume that the total space ${\mathcal  
D}=\cup_{t\in B}(t, D(t))$ is pseudoconvex in $B \times M$. If 
$\frac{\partial^2 \lambda }{ \partial t \partial \overline{ t} 
}(0)=0 $, then $\frac{\partial g (t,z)}{\partial t}\bigl|_{t=0}\equiv 
0$ on $D(0)$.}

\vspace{2mm} 
Next let $\psi (t,z)$ be a $C^\infty-$defining function of ${\mathcal  
D}$ in a neighborhood of $\partial 
{\mathcal D}=\cup_{t\in B}(t, \partial D(t))$. Since $\partial 
D(t)$ is smooth, we have 
$$\bigl(\frac{\partial \psi}{\partial z_1}(t,z),...,\frac{\partial \psi}{\partial z_n}(t,z) \bigr)\ne (0,...,0)$$ for $(t,z)\in \partial 
{\mathcal  D}=\{\psi(t,z) =0\}$. We have a type of contrapositive of $(\star 1)$:

\vspace{3mm} 
    
$(\star 2)$ \ \ Assume that ${\mathcal  D}$ is pseudoconvex in $B 
\times M$. If there exists a point $z_0\in \partial D(0)$ with  
\begin{align}
 \label{contrapositive}
\frac{\partial \psi }{\partial t}(0, z_0) \ne 0,
\end{align}
then $\frac{\partial^2 (-\lambda)}{\partial t \partial \overline{ t}}(0)>0 $.

\vspace{3mm} 
\noindent {\bf  Proof of $(\star 2)$.} \ 
 We set $z_0=(z_1^0, \ldots , z_n^0);\,  z_k^0=x_k^0+iy_k^0$
 and $t=t_1+ it_2$; we may assume $\frac{\partial \psi}{\partial y_1}(0, z_0)\ne 0$. 
We write  $z'=(z_2, \ldots , z_n)$ and $z_0'=(z_2^0, \ldots ,z_n^0)$. 
In a sufficiently small neighborhood $B_0 \times V =:B_0 \times (V_1 
\times V')$ of $(0, z_0)=(0, z_1^0, z'_0)$ we 
can write  $\partial D(t)$ in the form
$$
y_1= y_1(t, x_1, z'):= c_0(t,z') + c_1(t,z')(x_1-x_1^0)+ c_2(t,z')
({  x_1-x_1^0})^2 + \ldots 
$$
where $c_0(0, z_0')=y_1^0$. 
Using (\ref{contrapositive}) we may assume $\frac{\partial \psi}{\partial 
t_1}\ne 0 $ in $B_0 \times (V_1 \times V')$. 
By taking a smaller product set $B_0 \times V'$ if necessary we can also assume that $\frac{\partial c_0(t, z')}{\partial t_1} \ne 0 $ on $B_0 \times V'$. We set 
$A=\frac{\partial c_0 }{\partial t_1}(0, z_0')\ne 0 $. It follows 
that 
\begin{align*}
& y_1(t_1,x_1, z')-y_1(0,x_1, z')\\[2mm]
 &=
( c_0(t_1,z')-c_0(0,z')) +(c_1(t_1,z')-c_1(0,z')) (x_1-x_0)+ 
(x_1-x_1^0)^2+ \cdots \\[2mm]
&= t_1\ \biggl([ A + O(t_1,z')]+ [A_1 + O(t_1, z')](x_1-x_1^0)
+\ldots \biggr).
\end{align*}
We can find a  small interval $I_0:=[-r, r]$ on the 
$t_1$-axis;  a small polydisk $V_0'$ of center $z_0'$ in $\mathbb{ C}^{n-1}$, 
and a sufficiently small interval $J_0=[x_1^0- r_0, x_1^0+r_0]$ on the $x_1$-axis such that  
$$
|y_1(t_1,x_1,z')-y_1(0, x_1,z')| \ge \frac{ A}2\ |t_1|  \quad \mbox{ on $J_0 
\times (I_0 \times V_0')$}.
$$
Using this estimate, it follows from the boundary behavior of the $c$-Green function $g(t,z)$ 
and standard potential-theoretic arguments that 
$
\frac{\partial g(t,z)}{\partial t_1}\bigl|_{t=0}\not\equiv  0, 
$
and hence $\frac{\partial^2 (-\lambda)}{\partial t \partial \overline{ t}}(0)>0 $.
 \hfill $\Box$ 

\vspace{2mm}

Returning to the case of a pseudoconvex domain $D$ in ${\mathcal  H}$ with $C^\omega-$smooth boundary, we proved in Lemma \ref{lem:basic} that under certain hypotheses on $\partial D$ the function $-\lambda[z,w]$ is a plurisubharmonic exhaustion function for $D$. The next lemma shows that if $\partial D$ hits, but does not contain, one of the tori ${\bf  T}_a$ or ${\bf  T}_b$, and $D$ does not contain the other one, then $D$ is Stein.  

\begin{lemma}
 \label{last-lem} 
\ \ Let  $D$ be a pseudoconvex domain in ${\mathcal  H}$ with $C^\omega-$smooth boundary. If 
 $\emptyset \not =\partial D\cap {\bf  T}_a \not = {\bf  T}_a$ and 
$D \not \supset  {\bf  T}_b$, then $D$ is Stein (and similarly if ${\bf  T}_a$ and ${\bf  T}_b$ are switched).
\end{lemma}

\vspace{1mm} 

The condition $D \not \supset {\bf  T}_b$ separates  
into the following three cases:
$$
 (c1) \ \ \partial D \cap {\bf  T}_b= \emptyset, \quad ({c2}) \ \ \emptyset \ne \partial D \cap {\bf  T}_b \ne   {\bf  T}_b  \quad \mbox{ or } 
\quad   ({c3}) \ \ \partial D \cap {\bf  T}_b ={\bf  T}_b.
$$

\vspace{1mm} 
\noindent   {\bf Proof.}\ \  We first want to show that if $-\lambda[z,w]$ is not strictly plurisubharmonic in $D$, then there is point $p_0= [z_0, w_0]$ in $D^*$ at which $-\lambda [z,w]$ is not strictly plurisubharmonic; then we show this cannot occur so that $D$ is Stein. Let $\psi[z,w]$ be a 
defining function for $D$ defined in a neighborhood of $\partial D$. We divide the proof of the lemma in five cases related to $\psi[z,w]$ and the subcases $(c1), (c2), (c3)$ of the condition $D \not \supset {\bf  T}_b$.
\vspace{3mm} 

{\it $1^{st}$ case.} \quad 
Assume there exists  
 $[z_0,0]\in \partial D\cap \bf{ T}_a$ with $z_0\ne 0$ such that neither 
$\frac{\partial \psi}{\partial z} $ nor $\frac{\partial \psi}{\partial w} 
$ vanishes at $(z_0,0)$ {
and assume case $(c1)$.}
\vspace{3mm}

Using $(\star 2)$, we first prove the following fact in 
this $1^{st}$ case. Assume $(1,0) \in D \cap {\bf  
T}_a$. Then $-\lambda [z,w]$ is strictly subharmonic at 
$[1,0]$ in the direction ${\bf  a}=(0,1)$, i.e., 
$$
\frac{\partial^2 (-\lambda)}{\partial \tau \partial \overline{ 
\tau}} [1,\tau]|_{\tau=0} >0.
$$
To see this, we take a small disk $\delta:=\{|\tau|< r\} \subset 
\mathbb{ C}_\tau$ and consider the variation of domains  
$$
\mathfrak D: \ \tau \in \delta \to D(\tau):=D[1,\tau] \subset \mathbb{C}_Z^* \times 
\mathbb{C}^* _W.
$$
Note that 
$$
D(\tau)= \left\{
\begin{array}{lll}
 \widetilde D^* \times (1, 1/\tau)& \mbox{ if}  &  \tau\in \delta 
  \setminus \{0\};\\[2mm]
\widetilde D_a \times \mathbb{ C}_W^* & \mbox{  if}  &  \tau=0
\end{array}
\right\}
$$

\vspace{1mm}\noindent  
(recall $D \cap {\bf  T}_a= [D_a, 0]$). We let $\lambda (\tau)=\lambda 
{ [1,\tau]}$ denote the $c$-Robin 
constant for $\bigl(D(\tau),  (1, 1)\bigr)$. We set $\mathfrak{ D}:=\cup_{\tau\in \delta } (\tau, D(\tau))$ 
and $\partial \mathfrak{ D}=\cup_{\tau\in \delta } (\tau, \partial 
D(\tau))$. For $\tau\in \delta \setminus \{0\}$, we consider the 
automorphism 
$$
F_\tau: \ (z,w) \in \mathbb{ C}_z ^* \times \mathbb{ C}_w^* \to  (Z,W)= 
(z, \frac{ w}\tau) \in \mathbb{ C}_Z^* \times \mathbb{ C}_W^*.
$$
From the definition of $D(\tau)$, we have $D(\tau)= F_\tau(\widetilde D^*) 
$. 
 We let $\psi (z,w)$ be a defining function for $ \partial D$ in ${\mathcal 
 H}$; to avoid notational issues we also regard $\psi(z,w)$ as a defining function of 
 $\partial \widetilde D$. For $\tau\in \delta \setminus  \{0\}$ we set
$$
\Phi(\tau,  (Z,W)):= \psi (Z, \tau W)
$$
which is a defining function  for $\partial \mathfrak{ D}|_{\delta \setminus 
\{0\}}$. 
Setting $\Phi [0, (Z,W)]:= \psi (Z,0)$, we see that  
$\Phi[\tau, (Z,W)]$ becomes a smooth defining function for the entire set $\partial \mathfrak{ D}$. We focus on the special point $(z_0,1)$ in $\partial D(0)$. Then 
\begin{align*}
 \nabla _{(Z,W)} \Phi \bigl|_{(0, (z_0, 1))}&=
\bigl( \frac{\partial\Phi }{\partial Z}, \frac{\partial \Phi }{\partial W}
\bigr)\bigl|_{(0, (z_0,1))}= \bigl( \frac{\partial\psi }{\partial z}, \frac{\partial \psi 
 }{\partial w} \ \tau
\bigr)\bigl|_{(0, (z_0,1))}\\
&= \bigl( \frac{\partial\psi }{\partial z}(z_0,0), 0\bigr) \ne (0,0) \quad 
 \mbox{ by the 
 condition of the $1^{st}$ step.}  
\end{align*}
Similarly, 
\begin{align*}
 \frac{\partial \Phi}{\partial \tau}\bigl|_{(0, (z_0, 1))}&=
 \frac{\partial\psi}{\partial w}\ W \bigl|_{(0,(z_0,1))}\\
&= \frac{\partial \psi }{\partial w}(z_0, 0) \ne 0 \quad 
\mbox{ by the 
 condition of the $1^{st}$ step.}  
\end{align*}
It follows from $(\star 2)$ that $\frac{\partial^2 (-\lambda)}{\partial \tau \partial \overline{ 
\tau}} [1,\tau]|_{\tau=0} >0 $, as desired.

We next prove that 
$-\lambda [z,w]$ in $D$ is strictly subharmonic at $[1,0]$ in  
{\it any} direction ${\bf  a}=(a_1,a_2)\in \mathbb{ C}^2 \setminus \{0\}$ with 
$\|{\bf  a}\|=1$ and $a_1\ne 0$,\,i.e., 
\begin{equation}\label{taumu}
\frac{\partial^2 (-\lambda)}{\partial \tau \partial \overline{ 
\tau}} [1+a_1 \tau, a_2\tau]\bigl|_{\tau=0}\, >0.
\end{equation}
We use the same notation $\tau$ and $\psi(z,w)$ as in the case ${\bf  a}=(1,0)$. 
We 
consider the variation of domains 
$$
\mathfrak{G}: \ \tau\in \delta  \to G(\tau):= D[1+a_1\tau, a_2\tau] 
\subset \mathbb{C}_Z^* \times \mathbb{C}_W^*. 
$$
Note that
\begin{align*}
 G(\tau)&= \left\{
\begin{array}{lll}
 \widetilde D^* \times (1/(1+a_1\tau), 1/(a_2\tau))& \mbox{ if}  &  \tau\in \delta 
  \setminus \{0\};\\[2mm]
\widetilde D_a \times \mathbb{ C}_W^* & \mbox{  if}  &  \tau=0
\end{array}
\right\} \quad &\mbox{ in case $a_2\ne 0$},\\[2mm]
G(\tau)&= [\widetilde D_a \times (1/(1+ \tau))] \times \mathbb{ C}_W^*  \quad \mbox{ if $
\, \tau \in \delta$} \quad &\mbox{ in case $a_2=0$}.
\end{align*}
We let $\mu(\tau):=\lambda [1+ a_1\tau, a_2\tau]$ denote the $c$-Robin 
constant for $\bigl(G(\tau), (1, 1)\bigr)$. Our claim (\ref{taumu}) is that
$\frac{\partial^2 (-\mu)}{\partial \tau \partial \overline{ \tau }}(0)>0 $. 

We set $\mathfrak{ G}:=
\cup_{\tau\in \delta } (\tau, G(\tau))$ 
and $\partial \mathfrak{ G}=\cup_{\tau\in \delta } (\tau, \partial 
G(\tau))$. 
Since $\frac{\partial \psi}{\partial z}(z_0,0) \ne 0$ {  and $a_1\ne 
0$},
we can find a point $W_0\in \mathbb{ C}_W^*$ such that  
\begin{align*}
 a_1 z_0 \frac{\partial \psi}{\partial z}(z_0,0)+ 
a_2W_0 \frac{\partial \psi}{\partial w}(z_0, 0) \ne 0.
\end{align*}
We note that $(z_0, W_0)\in \partial G(0)=(\partial \widetilde D_a) \times \mathbb{ C}_W^*$.
 We consider 
$$
\Psi (\tau, (Z,W)):= \psi ((1+a_1\tau)Z, a_2\tau W),
$$
which is defined in a sufficiently small polydisk ${\mathcal  V}:= \delta_1 \times (U_1 \times V_1) $ of 
center $(0, (z_0, W_0)) $ 
in $\delta \times \mathbb{ C}_Z^* \times 
 \mathbb{ C}_W^*$. This is a 
defining function for $\partial \mathfrak{ G}$ in ${\mathcal  V}$.
We have 
\begin{align*}
 \nabla _{(Z,W)} \Psi \bigl|_{(0, (z_0, W_0))}&=
 \bigl( \frac{\partial\psi }{\partial z}\cdot (1+a_1\tau), \frac{\partial \psi 
 }{\partial w} \cdot  a_2\tau
\bigr)\bigl|_{(0, (z_0,W_0))}\\
&= \bigl( \frac{\partial\psi }{\partial z}(z_0,0), 0\bigr) \ne (0,0); \\[3mm]
 \frac{\partial \Psi}{\partial \tau}\ \bigl|_{(0, (z_0, W_0))}&=
 \frac{\partial\psi}{\partial z} \cdot (a_1Z)+ \frac{\partial \psi}{\partial 
 w}\cdot ( a_2 W) \biggl]_{(0,(z_0,W_0))}\\
&=  a_1z_0\frac{\partial \psi }{\partial z}(z_0, 0)+ a_2W_0\frac{ \partial 
 \psi}{\partial w}(z_0,0) \ne 0.
\end{align*}
Using $(\star 2)$ we conclude that 
$\frac{\partial^2 
(-\mu)
}
{\partial \tau \partial 
\overline{ \tau}} (0)>0$ which proves our claim.

A similar argument shows that $-\lambda [z,w]$ in $D$ is strictly 
plurisubharmonic at any point $[z,0]\in D \cap {\bf  T}_a$. Hence, in case $(c1)$, we conclude that if $-\lambda [z,w]$ is not strictly plurisubharmonic in $D$, there exists a point $p'= [z', w']$  in $D^*$ at which $-\lambda [z,w]$ is 
not strictly plurisubharmonic. Now since $ \partial D\not \supset {\bf  T}_a$ and 
       $\partial D\not \supset {\bf  T}_b$, we are in case (2) (i) of Lemma \ref{essen-lemm}. Hence we have $\partial D \cap ({\bf  T}_a \cup {\bf  T}_b) 
       =\emptyset$. This contradicts $\partial D \cap {\bf  T}_a\not 
= \emptyset$; thus $D$ is Stein. \hfill $\Box$ 

\vspace{3mm} 

$2^{nd} \ case.$ \  Assume there exists  
 $[z_0,0]\in \partial D\cap \bf{ T}_a$ with $z_0\ne 0$ such that neither 
$\frac{\partial \psi}{\partial z} $ nor $\frac{\partial \psi}{\partial w} 
$ vanishes at $(z_0,0)$ and 
there exists  
 $[0,w_0]\in \partial D\cap \bf{ T}_b$ with $w_0\ne 0$ such that neither 
$\frac{\partial \psi}{\partial z} $ nor $\frac{\partial \psi}{\partial w} 
$ vanishes at $(0, w_0)$, and assume case $(c2)$. 

\vspace{3mm} 
Using the same argument as in the $1^{st}$ case we see that $-\lambda[z,w]$  
is strictly plurisubharmonic at any point 
$[0,w]\in D \cap {\bf  T}_b$ and at any point $[z,0] \in D \cap {\bf  T}_a$. Thus there again exists a point $p'= [z', w']$  in $D^*$ at which $-\lambda [z,w]$ is 
not strictly plurisubharmonic; and we similarly conclude that  $D$ 
is Stein. \hfill $\Box$ 

\vspace{3mm} 
 
$3^{rd} \ case.$ Assume there exists  
 $[z_0,0]\in \partial D\cap \bf{ T}_a$ with $z_0\ne 0$ such that neither 
$\frac{\partial \psi}{\partial z} $ nor $\frac{\partial \psi}{\partial w} 
$ vanishes at $(z_0,0)$ {
and assume case $(c3)$. } 
\vspace{3mm} 

{
Recall $\partial D 
\supset {\bf  T}_b$ holds in case $(c3)$. Here we need the function $U[z,w]$ on ${\mathcal  H}^*$ defined in section 2. Using 2 (b) of Lemma \ref{lem:basic}, i.e., for    
$[z_0, w_0]\in \partial D\setminus {\bf  T}_b$, 
$$
-\lambda [z,w] \to \infty \quad \hbox{ as}  \ [z,w]\in D \to [z_0,w_0],
$$ and property (1) of $U[z,w]$ we see that 
\begin{align}
 \label{szw}
s[z,w]:= {\rm max}\{ \lambda [z,w],  - U[z,w] \}  
\end{align}
is a well-defined plurisubharmonic exhaustion function for $D$. 
In order to prove that $D$ is Stein, we use a result from $\S$ 14 in \cite{oka}: it suffices to show that  
for any $K \Subset D$ there exists a Stein domain $D_K$ with $K 
\Subset D_K \subset D$. To construct $D_K$, we take $m > \max_{[z,w]\in K}|-\lambda [z,w]|$ and 
consider 
\begin{align}
  \label{vzw}
v[z,w]:= {\rm max} \{ -\lambda [z,w]+2m, \   - \varepsilon U[z,w] \}
\end{align}
where $\varepsilon >0$ is chosen sufficiently small so that $v[z,w]= -\lambda [z,w]+2m $ on 
$K$. Again from property (1) of $U[z,w]$, $v[z,w]$ is a well-defined plurisubharmonic exhaustion function for $D$.
We take $M>1$ sufficiently large so that 
$$
K\Subset D(M):=\{[z,w]\in D: v[z,w] <M\} \ \hbox{and} \  \emptyset \ne \partial D(M)\cap {\bf  T}_a\not = {\bf  T}_a.
$$
Note that $D(M) \Subset D$; thus $\partial D \supset {\bf T}_b$ implies that ${\bf T}_b \cap \overline{D(M)}=\emptyset$; also $\partial D(M)$ is piecewise smooth. We now have 
\begin{equation}
\label{dmprop} 
\partial D(M)\cap {\bf  T}_b=\emptyset \ \hbox{and}  \  \emptyset \ne 
\partial D(M)\cap {\bf  T}_a\not = {\bf  T}_a.
\end{equation} 
We consider the $c$-Robin function $\lambda _M[z,w]$ for $D(M)$. Although 
$\partial D(M)$ is not smooth, by the construction of $\lambda 
_M[z,w]$ and the fact that $\partial D(M) \not \supset {\bf  T}_a, {\bf  T}_b $,  
it follows that $-\lambda_M[z,w]$ is a smooth plurisubharmonic exhaustion function for $D(M)$. 

Let $D(M,M'):=\{[z,w]\in D(M): -\lambda _M[z,w] <M'\}$ and take $M'>1$ sufficiently large so that   
$$
D(M,M') \Supset K \hbox{and}  \ \emptyset \ne \partial D(M,M')\cap  {\bf  T}_a \not \supset {\bf  T}_a. 
$$
Now since $-\lambda_M[z,w]$ is smooth we have that $D(M,M')$ is a pseudoconvex domain in ${\mathcal  H}$  
with smooth boundary; moreover we have
\begin{equation}\label{dmmprop}
\partial D(M,M')\cap
{\bf   T}_b=\emptyset \ \hbox{and}  \ \emptyset \ne \partial D(M,M')\cap  {\bf  T}_a \not \supset {\bf  T}_a. 
\end{equation} 
We can now apply the $1^{st}$ case, where we assumed condition $(c1)$, to $D(M,M')$ to conclude that $D(M,M')$ is Stein; hence $D$ is Stein.}

\vspace{3mm}

{\it $4^{th}$ case.} \  Assume one of $ \frac{\partial  \psi} {\partial z}, 
\frac{\partial  \psi }{\partial w}  $ vanishes identically 
on $\partial D\cap \bf{ T}_a$ 
 and assume case $(c1)$. 
\vspace{3mm} 

To deal with this case we  construct the 
  $C$-Robin function $\Lambda [z,w]$ on $D$ with respect to a positive 
  constant function $C$ on $\mathbb{ P}^2\supset \mathbb{C}^2$ and 
the restriction of the Fubini-Study metric  $dS^2$ on $\mathbb{ P}^2$ to $\mathbb{ C}^* \times \mathbb{ C}^*$. Note this metric is different than the Euclidean 
metric $ds^2$ on $\mathbb{ C}^2$ restricted to  
$\mathbb{ C}^*\times \mathbb{ C}^*$; accordingly, $-\Lambda[z,w]$ is a 
smooth plurisubharmonic exhaustion function on $D$ which is different from the function 
$-\lambda [z,w]$. Moreover, for any positive constant $k$ the function $u_k[z,w]:= -(\lambda [z,w]+ k \Lambda [z,w]) 
$ is a smooth plurisubharmonic exhaustion function for $D$. We claim that we can find a $k$ and an increasing sequence 
$\{M_n\}_{n=1,2,\ldots}$ tending to $+\infty$ such that the increasing sequence of pseudoconvex domains $D_n= \{ 
[z,w]\in D:
u_k[z,w]<M_n\}$ satisfy the hypotheses of the $1^{st}$ case. Clearly $\partial D_n\cap \bf{ T}_b=\emptyset$ so that $(c1)$ holds. It remains to select $k$ and then the sequence $M_n$ so that there exists  
 $[z_n,0]\in \partial D\cap \bf{ T}_a$ with $z_n\ne 0$ such that neither 
$\frac{\partial \psi_n}{\partial z} $ nor $\frac{\partial \psi_n}{\partial w} 
$ vanishes at $(z_n,0)$ where $\psi_n[z,w]:= u_k[z,w]-M_n$. From the $1^{st}$ case we conclude that each $D_n$ is Stein and it follows from $\S$ 14 of \cite{oka} that $D$ is Stein.

\vspace{3mm} 
  {\it $5^{th}$ case.} \  Assume one of $ \frac{\partial  \psi} {\partial z}, 
\frac{\partial  \psi }{\partial w}  $ vanishes identically 
on $\partial D\cap \bf{ T}_a$ 
 and assume case $(c2)$ or $(c3)$.

\vspace{1mm} 
The type of argument used to show a domain $D$ in the 
$2^{nd}$ or $3^{rd}$ case, where we assume $(c2)$ or $(c3)$ of the condition $D \not \supset {\bf  T}_b$, reduces to the $1^{st}$ case, where we assume $(c1)$ of this condition, allows us to deduce the $5^{th}$ case from the $4^{th}$ case. We leave the details to the reader. 
\hfill $\Box$ 

\vspace{1mm}

We next turn to the situation where $\partial D$ contains one of ${\bf  T}_a$ or ${\bf  T}_b$ but not both.
\begin{lemma} \label{nemtype}
\ \ Let  $D$ be a pseudoconvex domain in ${\mathcal  H}$ with $C^\omega-$smooth boundary. If 
$ (i) \ \partial D \supset {\bf  T}_a$ and $(ii)\  \partial D 
\cap  {\bf  T}_b  \not =  {\bf  T}_b$, then
\begin{enumerate}
 \item [(1)]  $D$ is  Stein or
 \item [(2)]  $D$ is  of 
 type {\rm (b2)} in Theorem \ref{thm-1}. In fact, $D=\bigcup_{c\in \delta} \sigma_c$ 
       with  $0\in \partial \delta$ and $\infty \not \in \delta\cup \partial \delta 
       $.
\end{enumerate} 
(and similarly if ${\bf  T}_a$ and ${\bf  T}_b$ are switched as well as $0$ and $\infty$). 
\end{lemma}

The condition (ii) separates into the following two cases:
$$
(\tilde c\, 1)  \ \ \emptyset \ne \partial D \cap {\bf  T}_b \ne {\bf  T}_b \quad 
\mbox{ or } \quad 
(\tilde c\, 2) \ \  D \supset {\bf  T}_b.
$$

\noindent {\bf Proof.} {
We first treat the case $(\tilde c\,1)$.} We assume that  $D$ is not of 
 type ({\rm b2}) as in (2) and we show $D$ is Stein. We proceed as in the proof of the $3^{rd}$ case of Lemma \ref{last-lem} where we use the function $U[z,w]$ on ${\mathcal  H}^*$ defined in section 2. However, instead of (\ref{szw}) and (\ref{vzw}) we use 
\begin{align*}
s[z,w]:= {\rm max}\{ -\lambda [z,w],   U[z,w] \}.  
\end{align*}
and 
$$
v[z,w]:= {\rm max} \{ -\lambda [z,w]+2m, \  \varepsilon U[z,w] \}
$$
We leave the details to the reader.
 
We next treat the case $(\tilde c\, 2)$ in which $\partial D 
\supset {\bf  T}_a$ and $D \supset {\bf  T}_b$. In this setting we shall show that conclusion (2) 
in Lemma \ref{nemtype} holds. 

Since ${\bf  T}_b$ is compact in $D$, we can find a neighborhood $V$ of ${\bf  T}_b$ in $D$ such that  
${\bf  T}_b \Subset V \Subset D$. Since $\Sigma _c:=\{|w|= 
c|z|^\rho\}{ / \sim}$ (or $\sigma_c:=\{ w= c z^\rho\}{  /\sim}$) approaches ${\bf  T}_b$ in 
${\mathcal  H}$ as $c \to \infty$, it follows that for $c$ sufficiently large, the Levi flat hypersurface $ \Sigma_c$ satisfies $\Sigma_c \Subset V \Subset D$ (or the compact torus $\sigma_c$ satisfies $\sigma_c\Subset V \Subset D$) . But $-\lambda [z,w]$ is a 
plurisubharmonic function on $D$ (although not necessarily an exhaustion function); hence $-\lambda [z,w]$ is not strictly plurisubharmonic at any point in $\Sigma _c$ (or $\sigma_c$). From Lemma \ref{essen-lemm}, we conclude that $D$ is given as in case (2) (ii) of that lemma. 

For simplicity, we complete the argument if $\Sigma_c \Subset V \Subset D$. 
We claim that $\rho$ is of case ({\rm b 2}) ($\rho$ rational and $\tau$ rational) in Theorem \ref{thm-1} and hence $D$ is of the form in  case (2) (ii-a) of Lemma \ref{essen-lemm}, completing our proof.
 For if $\rho$ is of case 
${\bf  a}$ ($\rho$ irrational) or of case ({\rm  b1}) ($\rho$ rational and $\tau$ irrational), then from the proof of Lemma \ref{essen-lemm}, we have (recall $(\alpha^* )$)
$$D^* =\bigcup _{c\in I}\, \Sigma_c =\bigcup _{c\in I}\, \{|w|=c |z|^\rho\}$$ where 
$I=(r,R)$ is an open interval in $(0, \infty)$ because $D^*$ is connected. Since $D\supset {\bf T}_b$, 
$ D= \cup _{c\in (r,\infty]} \Sigma _c$. However, since $\partial D \supset {\bf  T}_a$, we must have $r=0$. Thus $D={\mathcal H}\setminus {\bf  T}_a$ which contradicts the smoothness of $\partial D$.   \hfill $\Box$

\vspace{3mm} 

Note in particular we have proved that the Nemirovskii-type domains in (2) (ii-b) of 
Lemma \ref{essen-lemm} are Stein. An entirely similar proof, which we omit, deals with the case where $\partial D$ contains both ${\bf  T}_a$ and ${\bf  T}_b$.

\begin{lemma} \label{finlem}
\  Let  $D$ be a pseudoconvex domain in ${\mathcal  H}$ with $C^\omega-$smooth boundary. If 
$ \partial D \supset {\bf  T}_a\cup {\bf T}_b$, then
\begin{enumerate}
 \item [(1)]  $D$ is  Stein or
 \item [(2)]  $D$ is  of 
 type {\rm (b2)} in Theorem \ref{thm-1}. In fact, $D=\bigcup_{c\in \delta} \sigma_c$ 
       with $0,\infty \in \partial \delta$.
\end{enumerate} 
\end{lemma}

\noindent We suspect that under the hypotheses of Lemma \ref{finlem} conclusion (2) must always hold, but we are unable to verify this.

We can now easily conclude with the proof of our main result. 

\vspace{3mm} 

\noindent {\bf Proof of Theorem \ref{thm-1}.} 
Let $D$ be a pseudoconvex domain in ${\mathcal  H}$ with $C^\omega-$smooth 
boundary which is not Stein. {
We consider three ``symmetric'' cases depending on the nature of $\partial D\cap {\bf  T}_a$ or $\partial D\cap {\bf  T}_b$. {

\vspace{3mm}

{\it $1^{st}$ case:} $\partial D\supset {\bf  T}_a$ (or $\partial D\supset {\bf  T}_b$).
\vspace{3mm} 

If $\partial D\supset {\bf  T}_a$, we can have either $\partial D 
\cap  {\bf  T}_b  \not =  {\bf  T}_b$ or $\partial D 
\supset  {\bf  T}_b$. If $\partial D 
\cap  {\bf  T}_b  \not =  {\bf  T}_b$, from Lemma \ref{nemtype}, $D=\bigcup_{c\in \delta} \sigma_c$ 
       with $0\in \partial \delta$ and $\infty \not \in \partial \delta $. If $\partial D 
\supset  {\bf  T}_b$, this means $ \partial D \supset {\bf  T}_a\cup {\bf T}_b$; hence Lemma \ref{finlem} implies $D=\bigcup_{c\in \delta} \sigma_c$ 
       with $0,\infty \in \partial \delta$. 
       
  \vspace{3mm} 
{\it $2^{nd}$ case:} $\partial D\cap {\bf  T}_a=\emptyset$ (or $\partial D\cap {\bf  T}_b=\emptyset$).

\vspace{3mm} 

If $\partial D\cap {\bf  T}_a=\emptyset$, we can have either $\partial D 
\cap  {\bf  T}_b  \not =  {\bf  T}_b$ or $\partial D 
\supset  {\bf  T}_b$. If $\partial D 
\supset  {\bf  T}_b$, we are done by the $1^{st}$ case. If $\partial D 
\cap  {\bf  T}_b  \not =  {\bf  T}_b$, either 
$$ (I)  \ \partial D 
\cap  {\bf  T}_b  =\emptyset \ \hbox{or} \   (II) \ \emptyset \not = \partial D\cap {\bf  T}_b\not ={\bf  T}_b.$$ 
Note that if $\partial D \cap {\bf  T}_b=\emptyset$, then in this $2^{nd}$ case 
$\partial D \cap ({\bf  T}_a \cup {\bf  T}_b)=\emptyset$.

Let $\lambda [z,w]$ be the $c$-Robin function of $D$. 
From Lemma \ref{lem:basic} we know that $-\lambda [z,w]$ is a plurisubharmonic 
exhaustion function on $D$. 
We shall prove that under our assumption that $D$ is not Stein we can find a point $[z_0,w_0]$ in $D^*$ at which $-\lambda [z,w]$ is not strictly plurisubharmonic. We give the proof when $\rho$ is irrational since the other cases are completely analogous.

In the setting of the $2^{nd}$ case with $ (I)  \ \partial D 
\cap  {\bf  T}_b  =\emptyset$ we have three possible situations for $D$ relative to ${\bf  T}_a, {\bf  T}_b$: $(i) \ D\cap ({\bf  T}_a  
\cup {\bf  T}_b)=\emptyset;$  
$(ii) \ D \cap {\bf  T}_a=\emptyset  \ \mbox{ and }  D \supset {\bf  T}_b \ 
\mbox{(or the symmetric case with ${\bf  T}_a, {\bf  T}_b$ switched)}$;  and 
$(iii) \ D \supset {\bf  
T}_a \cup {\bf  T}_b.$ 

In case (i) we are done since $D=D^*$ so that, by the assumption $D$ is not Stein, there is a point $[z_0,w_0]$ in $D=D^*$ at which $-\lambda [z,w]$ is not strictly plurisubharmonic. By (2) (i) of Lemma \ref{essen-lemm}, $D$ is of type 
{\rm (a-1)}. In case (ii), since ${\bf  T}_b$ is compact in $D$, we can find a neighborhood $V$ of ${\bf  T}_b$ in $D$ such that  
${\bf  T}_b \Subset V \Subset D$. The Levi flat hypersurface $\Sigma _c:=\{|w|= c|z|^\rho\}$ approaches ${\bf  T}_b$ as $c \to \infty$; hence $ \Sigma_c 
\Subset V \Subset D$ for $c$ sufficiently large. Since $-\lambda [z,w]$ is a 
plurisubharmonic function on $D$, $-\lambda [z,w]$ is not strictly 
plurisubharmonic at points of $\Sigma _c$; thus we can find such a point in $D^*$. Recalling $(\alpha^*)$:
$$D^* =\bigcup _{c\in I}\ \{|w|=c |z|^\rho\},  \mbox{ where 
$I$ is an open interval in} \ (0, \infty),$$
we see that $D$ is of type (a-2'') in Theorem \ref{thm-1}. 
In case (iii), similar reasoning as in case (ii) shows that $\Sigma _{c_0}\subset D $ for some $c_0\ne 0, \infty$. It follows 
 that $D= \bigcup_{c\in I} \Sigma_c$ where $I $ is an interval in $[0, 
 \infty]$. Since $D \supset {\bf  T}_a\cup {\bf  T}_b$, we have $I=[0, \infty]$, 
 i.e., $D= {\mathcal  H}$, which is absurd. This finishes the proof of of the $2^{nd}$ case under situation $(I)$.

\vspace{3mm} 

{
To finish the proof of the $2^{nd}$ case, where $\partial D \cap {\bf  T}_a =\emptyset$, it remains to deal with situation $(II)$, i.e.,  $\partial D\cap {\bf  T}_a=\emptyset$ and 
$\emptyset \not = \partial D\cap {\bf  T}_b\not ={\bf  T}_b$. Apriori, we 
separate this into two subcases: 
 $$(c1)\  D \supset {\bf  T}_a \ \hbox{and} \  (c2) \ D \not \supset {\bf  T}_a.$$ 
 In case (c1), using the argument in case 
 (ii) above we can find a neighborhood $V$ of ${\bf  T}_a$ in $D$ such that ${\bf  T}_b \Subset V \Subset D$ and hence $ \Sigma_c 
\Subset V \Subset D$ for $c>0$ sufficiently close to $0$. Thus we obtain points in $D^*$ at which $-\lambda [z,w]$ is not strictly 
plurisubharmonic. We now appeal to case (2) (i) of Lemma \ref{essen-lemm}.   

Now we observe that case (c2) cannot occur, for the assumptions $\emptyset \not = \partial D\cap {\bf  T}_b\not ={\bf  T}_b$ and $D \not \supset {\bf  T}_a$ imply from Lemma \ref{last-lem} that $D$ is Stein.
  }

     \vspace{3mm} 
{\it $3^{rd}$ case:} $\emptyset \not = \partial D\cap {\bf  T}_a\not ={\bf  T}_a$ (or $\emptyset \not = \partial D\cap {\bf  T}_b\not ={\bf  T}_b$).

\vspace{3mm} 

If $\emptyset \not = \partial D\cap {\bf  T}_a\not ={\bf  T}_a$, from Lemma \ref{last-lem} we must have  $D\supset {\bf T}_b$. Thus $\partial D\cap {\bf  T}_b=\emptyset$ and we are done by the $2^{nd}$ case. 

This completes the proof of Theorem \ref{thm-1}. \hfill $\Box$ 

 \vspace{3mm} We end with an explicit example of the construction of both $D[z,w]$ and the $c-$Robin function $\lambda [z,w]$ for a specific Nemirovskii-type domain $D\subset \mathcal H$. We recall the fundamental domain ${\mathcal  
F}= E_1\cup E_2=(E_1'\cup E_1'')\cup(E_2'\cup E_2'')$ for $\mathcal H$ defined in (\ref{eqn:found}). Let $D$ be a 
subdomain of ${\mathcal  F}$ defined by 
$$D:=(E_1' \times K_1'')\cup (E_2' \times K_2'') \subset E_1\cup E_2$$
where (recall $b>1$)
$$K_1'':=\{1<|w| < b\}\cap \{\Re w >0\} \ \hbox{and} \ K_2'':=\{|w| < b\}\cap \{\Re w >0\}.$$
We note that $\partial D$, which can be written as
$$ \{|z|\le |a|\} \times \{\Re\, w=0, 1\le|w|\le b\} 
 \bigcup \{1\le |z| \le |a|\} \times 
\{\Re\, w=0, |w|\le 1\},$$
is smooth in ${\mathcal  H}$. To see that $D$ is of Nemirovskii-type as in Lemma \ref{essen-lemm} (ii-b),  
setting $$N=\mathbb{ C}_z \times \{\Re\, w>0\} \subset (\mathbb{ 
C}^2)^*$$ we will show that  
\begin{align}
 \label{nemi-d}
N/ \sim  \ = D \  \mbox{ in ${\mathcal  H}$},\ \ \ \mbox{ or equivalently},  
 \ \ \
N= \widetilde D  = D \times {\mathcal  I} \mbox{ in $(\mathbb{ C}^2)^*$}
\end{align}
(recall (\ref{ai})). Hence  $N \setminus \mathbb{ C}_z \times \{0\}= \widetilde {D^*}$. 

\vspace{3mm} 
To prove (\ref{nemi-d}), we show $N= \widetilde D$.  
Let $(z,w)\in N$. Then we have $z= a^nz_0$ and $w=b^mw_0$ for some 
$n,m\in \mathbb{ Z}$ and $(z_0,w_0)\in {\mathcal  F}$. Since $b>1$, we 
have $\Re\, {w_0}>0$. 
\vspace{2mm}

 {\it Case 1: $\ n \ge m$.}
\vspace{2mm}

In this case we have $(z,w) \sim (z/a^n, w/b^n)= (z_0, b^{m-n}w_0) \in E_2' \times K_2'' \subset D$.


\vspace{2mm}
{\it Case 2: $m \ge n$.}
\vspace{2mm}

In this case we have $(z,w) \sim (z/a^m, w/b^m)= (a^{n-m}z_0, w_0) \in E_1' 
\times K_1'' \subset D$.
\vspace{2mm}

\noindent Hence $N \subset \widetilde D=D \times {\mathcal  I}$. The converse is 
clear from the relations $D \subset N$ and $N {\mathcal  I}= N$.  

\vspace{3mm} 

We turn to the study of the sets $D[z,w]$ and the $c-$Robin functions $\lambda [z,w]$ for 
$(D[z,w], e)$ with respect to the metric $ds^2$ on $\mathbb{C}^* \times 
\mathbb{C}^*$ and the function $c(z,w)>0$. Recall $e=(1,1)$. Note that $\widetilde {K_1''}= \{ \Re\, w >0\}$. Let $w' \in K_2''$. We write $w'= |w'| e ^{i\theta}$ where $-\frac{ \pi}2< \theta < \frac{ \pi}2$ and define
\begin{equation} \label{Delta} \delta (w'):= \{ w=u+iv\in \mathbb{ C}_w : (\cos \theta )u- (\sin \theta) v>0\}.
\end{equation}
We then have 
$$
\{ \Re\, w >0 \} \times \frac{ 1}{w'}= \delta (w') \quad \mbox{ in $\mathbb{ C}_w$},
$$
so that dist$(1,\partial \delta (w'))\geq \cos \theta$ for $|w'|\leq 1$. Recalling the formulas
\begin{equation*}
 \begin{array}{lll}
  & D[z,w]= \bigl(\displaystyle{ (\frac{ 1}z, \frac{ 1}w)\cdot D}\bigr) \times {\mathcal 
 I} & \ \  \mbox{if } \ [z,w] \in D^*;\\[2mm]
&D[z,0]= \bigl(\displaystyle{ \frac{ 1}z D_a, \mathbb{ C}^*}\bigr) \times {\mathcal  I}
  \ =(\displaystyle{ \frac{ 1}z}\, \widetilde { D_a}) \times \mathbb{ C}_w^* &  \ \ \mbox{if } \ [z,0] \in  D\cap {\bf T}_a;\\ [2mm]
&D[0,w]= \bigl(\displaystyle{ \mathbb{ C}^*, \frac{ 1}w D_b} \bigr) \times {\mathcal  I}
 \ = \mathbb{ C}_z^* 
\times \displaystyle{ (\frac{ 1}w}\, \widetilde {D_b})  
& \ \ \mbox{if } \ [0,w] \in D\cap  {\bf  T}_b 
\end{array}\end{equation*}
where $D \cap {\bf  T}_a= D_a \times 
\{0\}$, $D \cap {\bf  T}_b= \{0\} \times D_b$, 
$\widetilde { D_a}=\{ a^nz: z\in D_a, \ n\in \mathbb{Z}\} \subset \mathbb{ C}_z^*$ and 
$\widetilde {D_b}=\{ b^nw : w\in D_b, \ n\in \mathbb{Z}\} \subset \mathbb{ C}_w^*$, in using the equality $\widetilde D=N$ we obtain the following:
\vspace{2mm}

If $(z',w') \in D^*$, then 
$$
D[z',w'] = (\frac{ 1}{z'},\frac{ 1}{w'})\, \widetilde{ D^*}= \mathbb{ C}_z^* \times \delta (w'),
$$
while if $(0, w')\in D$, then 
$$
D[0,w']= \mathbb{ C}_z^* \times \frac{ 1}{w'}\,  \widetilde {K_1''} =
\mathbb{ C}_z^* \times \delta (w'). 
$$
{\it Hence for any $[z,w] \in D$, we have
$$
D[z,w]= \mathbb{ C}_z^* \times \delta (w)
$$
which is independent of $z$. It follows that $\lambda [z,w], \ [z,w]\in D$ is independent of $z$.}   

\vspace{2mm}

We analyze the boundary behavior of $\lambda [z,w]$. We consider different cases:
\begin{enumerate}
 \item [(1)] Let $ [z_0, w_0]\in \partial 
D \setminus {\bf T}_a$; i.e., $z_0 \ne 0, \ w_0= 0+ i v_0 \ne 0$. 
We let $[z,w]\in D$ approach $[z_0, iv_0]$. If $z\to z_0$ 
and $w \to iv_0$, then regarding (\ref{Delta}) with  
       $\theta= \pi/2$ we see that 
\begin{align*}
 D[z,w] =\mathbb{ C}_z^* \times \delta (w) \ \hbox{approaches} \ D[z_0, iv_0] =\mathbb{ C}_z^* \times 
\{\Im\, w<0\}.
\end{align*}
In particular $e\in 
\partial (\mathbb{ C}_z^* \times \{\Im\, w<0\})$; thus as $[z,w]$ approaches $[z_0, iv_0]$, we have 
${\rm dist} \bigl(\partial D[z,w], e\bigr)$ tends to $0$ and $\lambda [z,w]$ tends to 
$-\infty$. 
 \item [(2)] Let $[z_0, 0] \in \partial D \cap {\bf  T}_a= {\bf  
       T}_a$ where $z_0 \ne 0$.  We let $[z,w]\in D$ approach $[z_0, 0] $ in such a way 
 that $z\to z_0$ arbitrarily but $w \to 0$ in an angular sector; i.e., writing $w=|w|e^{i\theta}$, there exists $\theta_0$ with $0<\theta_0< \pi/2$ so that $|\theta|<\theta_0$ as $|w|\to 0$. As before we have $D[z,w]= \mathbb{ C}_z^* \times \delta (w)$. It follows from (\ref{Delta}) that 
 dist$(\partial D[z,w],e)\geq \cos \theta_0$ for $|w|\leq 1$. Let $A$ be the $c-$Robin constant for 
the region $$G(\theta_0):=\{(z,w)\in \mathbb{ C}_z^* \times\mathbb{ C}_w^*: \  |z-1|^2+|w-1|^2< \cos^2\theta_0\}$$ with pole $e$. Then $A$ is finite and since $G(\theta_0)\subset D[z,w]$ for $|w|\leq 1$, clearly $\lambda[z,w]>A$. Thus $-\lambda [z,w]$ is not an exhaustion function due to its boundary behavior at ${\bf T}_a$. 
\end{enumerate}

Finally, we let $X:= w \frac{\partial }{\partial w} $ and $p_0=[z_0,w_0]\in D^*$.
Then the integral curve for $X$ with initial value $p_0$ is given by 
$$\sigma :=
p_0 \exp tX= (\{ \mathbb{ C}_z^*\} \times 
\{w_0\}) /\sim \ \subset \widetilde {D^*}/\sim \ = D^*.$$ Thus this example does indeed satisfy (1) and (2) (ii-b) of Lemma 4.2.

}

\section{Appendix A: Proof of Lemma \ref{fund-1}}

We give the proof of Lemma \ref{fund-1}. Assertion 1. follows from 
property (2) of the function $U[z,w]$ at the end of section 2. To see this, 
by definiton we note that $\widetilde \Sigma :=\{ U[z,w]=0\}$ coincides with $\{ |w|= 
|z|^\rho\}/\sim$ in ${\mathcal  H}^*$. We consider the integral curve $\widetilde \sigma_u$ 
of $X_u$ with initial value $ [1,1]$, i.e., $\widetilde \sigma_u 
=\{w= z^\rho\}/ \sim$. Since $\rho$ is real, we have 
 $ \widetilde \sigma_u \subset \{|w|= |z|^\rho\}/\sim \ =\widetilde 
 \Sigma$; hence $\widetilde \Sigma _u \subset \widetilde \Sigma$. 
Assume that $\rho$ is irrational, given $z_0\in \mathbb{ C}_z^*$ with 
$1<|z_0|\le |a|$, writing $pr\{z_0^{\rho}\}:= |z_0|^{\rho} e ^{i \rho \theta}$ for some $\theta$ we have 
$$
z_0^\rho =\{\,  pr\{z_0^\rho\}\, e^{2\pi i \rho n}: n\in \mathbb{ Z}\,\} \quad 
\mbox{ as sets in $\mathbb{ C}_w^*$}.
$$
This set is dense in the circle $\{|w|=|z_0|^\rho\}$. It follows that 
$\widetilde \sigma_u$ is dense in $\widetilde \Sigma $; hence $\widetilde \Sigma _u= 
\widetilde \Sigma $. This proves Assertion 1.(1) if $\rho$ is irrational, and it also proves Assertion $(\alpha)$ listed at the conclusion of Lemma \ref{fund-1} in this case. 

We next prove 1.(1) assuming $\tau$ is irrational and $\rho=q/p$ is rational. Again writing $pr\{z^{q/p}\}:= |z|^{q/p} e ^{i (q/p)\theta}$, we have 
\begin{align} \label{mendo}
\nonumber  \widetilde \sigma_u &=\{w=z^{q/p}\}/ \sim  \quad \mbox{ (by analytic 
 continuation)} \\ 
\nonumber &= \cup_{n \in \mathbb{ Z}} 
\{  (a^nz,(a^nz)^{q/p}: z \in \mathbb{ C}^*\} / \sim \\ 
\nonumber &= \cup_{n \in \mathbb{ Z}} 
\{  (z,  b^{-n} ((a^nz)^{q/p} ): z\in \mathbb{ C}^*\} / \sim \\ 
&=
 \cup_{n,k \in \mathbb{ Z}} \{(z, pr\{z^{q/p}\}\, e^{2\pi i 
 ((nq/p)+ k \tau)}:z\in \mathbb{ C}^*\}/\sim \ .
\end{align}
Since $\tau$ is irrational, we similarly  have 
 $\widetilde \Sigma_u = \widetilde \Sigma$, 
finishing the proof of 1.(1). A similar argument yields Assertion 
$(\alpha)$ in this case, completing its proof as well.

We next prove 1.(2) assuming $\tau:= ((q/p) \arg a- \arg b)/2\pi$ from (\ref{taueqn}) is rational (see Case $({\rm b2})$ in Theorem \ref{thm-1}). We defined  
$\mathfrak{ X}= \{\alpha z \frac{\partial }{\partial z}+ 
\beta { w} \frac{ \partial }{\partial w}: 
\alpha , \beta\in \mathbb{ C} \}$ on  ${\mathcal  H}^*$, 
which is a two-dimensional Lie algebra 
in ${\mathcal  H}^*$, and 
$\mathfrak{ X}_u=\{c X_u: c\in \mathbb{ C}\}$, 
which is a one-dimensional Lie subalgebra of $\mathfrak{ X}$. 
Then $\widetilde \sigma _u$ coincides with 
the Lie subgroup of ${\mathcal  H}^*$ corresponding to 
$\mathfrak{ X}_u$. We will give a concrete description of $\widetilde \sigma _u$ as a compact curve in ${\mathcal  H}^*$.

We let $g$ be the greatest common divisor of $p$ and $l$ and define $\nu:=pl/g\in \mathbb{ Z}$. Then 
$$
\{(nq/p) +k\tau:  n,k \in \mathbb{ Z}\}= \{j/\nu: 0 \le j \le \nu-1\} 
 \ \mbox{ mod $\mathbb{ Z}$}.
$$
It follows from (\ref{mendo}) that 
\begin{align*}
\widetilde \sigma _u&= \cup_{k=0}^{\nu-1} \{(z, pr\{z^{q/p}\}\, e^{2\pi i k/\nu})
:z\in \mathbb{ C}^*\}/\sim.
\end{align*}
Setting $
w_k(z): =pr\{z^{q/p}\}\, e^{2\pi i k/\nu} \in \mathbb{ C}_w^*, \  
k=0,1,\ldots ,\nu-1$ and ${\mathcal  W(z)}:= \cup_{k=0}^{\nu-1}  w_k(z) 
\subset \mathbb{ C}_w^*$, 
we have 
$$\widetilde \sigma _u=  \cup _{z \in \mathbb{ C}^*} (z, {\mathcal  W(z)})/\sim.$$
For $n\in 
\mathbb{ Z}$ we have $$(z, 
{\mathcal  W}(z))/\sim \ = (a^nz,  {\mathcal  W}(a^nz))/\sim.$$  
Moreover, since $1 <|w_k(z)|\le |b|$ for $1< |z|\le |a|$, 
$\widetilde 
\sigma _u$ may be considered as a graph in the fundamental domain ${\mathcal  F}$ 
(or as a multi-valued function $w={\mathcal  W}(z)$ on the annulus 
$\{1<|z|\le |a|$\} ) in the following manner:
$$
(*) \qquad \widetilde \sigma _u= \cup _{1<|z| \le |a|} (z, {\mathcal  W(z)}).  
$$
Since  $(z,{\mathcal  W}(z))\equiv (az, {\mathcal  W}(az))$ in ${\mathcal 
 H}$, $\widetilde \sigma _u$ is a compact curve in ${\mathcal  H}^*$; indeed, $\widetilde \sigma _u$ is a one-dimensional torus.

Moreover, if we consider the 
finite subgroup  $K= \{e^{2i\pi k/ \nu}: k=0, 1,\ldots ,  \nu-1\}$ in 
$\mathbb{ C}^*$, then for $c, c'\in \mathbb{ C}^*$ we have  
\begin{equation}\label{lastone}
\{w=c z^{q/p}\}/ \sim \ = \{ w= c'z^{q/p}\}/\sim 
 \quad \mbox{ if and only if \ $c'\in cK$}.
\end{equation}

To verify (\ref{lastone}), let $w_0=c z^{q/p}_0$ where 
$z_0\ne 0$. We can find $(z_0', w_0') \in \mathbb{ C}^* \times \mathbb{C}^*$
with $w_0'= c'(z_0') ^{q/p}$  and $(z_0', w_0')= (a^n z_0,b^nw_0)$ 
for some $n\in \mathbb{ Z}$. Then $cb^nz_0^{q/p}=c'(a^n z_0)^ {q/p}$; hence $|c|=|c'|$. Consider the 
 fundamental domain 
${\mathcal  F}_{|c|}:={\mathcal  F}\times (1, |c|)$ of ${\mathcal  H}$. 
Similar to $(*)$, we have the following equalities in ${\mathcal  
F}_{|c|}$: 
\begin{align*}
 \{w= c z^{q/p}\}/\sim &= \cup _{1<|z| \le |a|} (z, c\, 
 {\mathcal  W}(z));  
\\
\{w=c' z^{q/p}\}/\sim &= \cup _{1< |z| \le |a|} (z, c'\, 
 {\mathcal  W}(z)).  
\end{align*}
Thus $\{w= c z^{q/p}\}/\sim \ =  \{w= c' z^{q/p}\}/\sim$ iff 
$c\, {\mathcal  W}(z)= c'\,{\mathcal  W}(z)$ as sets in $\mathbb{ C}^*_w$ for any 
$z\in \{1<|z|\le |a|\}$. This is clearly equivalent to $c'\in cK$.

Thus we can write ${\mathcal  H}^*$ as a disjoint union:
$${\mathcal  H}^*= \bigcup_{cK \in \mathbb{  C}^*/K} \{w=c z^{q/p}\}/ \sim. 
$$
Since $\widetilde \sigma _u$ is compact, we have 
$$
\lim_{ c\to 0} \ \{w=c z^{q/p}\}/ \sim\ = \ {\bf  T}_a \quad \hbox{and}  \quad  
\lim_{ c\to \infty} \ \{w=c z^{q/p}\}/ \sim\ =\ {\bf  T}_b
$$
in ${\mathcal  H}$. Now since $\mathbb{  C}^*/K$ is equivalent to $\mathbb{ 
C}^*$, we write $\mathbb{  C}^*/K= \mathbb{ C}^*$; $cK= c$; 
$\{w=c z^{q/p}\}/ \sim \ = \sigma_c$; ${\bf  T}_a= \sigma _0$, and ${\bf  
T}_b= \sigma _\infty$. With this notation can write ${\mathcal  H}$ as a disjoint union 
$$ {\mathcal  H}= \cup _{c\in 
\mathbb{ P}^1} \ \sigma _c.$$ 
This proves Assertion 1.(2) and $(\beta)$.

 We now prove 2.(1).  
Let  $X=\alpha z \frac{\partial }{\partial z}+ \beta w \frac{\partial 
 }{\partial w}  \not \in \{c X_u: c\in \mathbb{ C}\}$ with $\alpha ,\beta 
 \ne 0$. Considering $X$ as a vector field in $\mathbb{ C}_z^* \times 
\mathbb{C}_w^*$, the integral curve $\{\exp tX: t\in \mathbb{ 
C}\}$  of $X$ with initial value $e=(1,1)$ in $\mathbb{ C}_z ^* \times \mathbb{ C}_w^*$ is 
$w= z^{\beta /\alpha } $.
Let $\beta /\alpha = A+Bi$ where $A,\,B$ are real. 
Then 
$$
w=z^{A+Bi}= e^{(A+Bi)\log z}.
$$
Fix $z'\in \mathbb{ C}^*$ and let ${  {\rm Log}\,} z'= \log |z'|+ i \theta'$ 
($0\le \theta'< 2\pi$) be the principal value. By analytic 
continuation, over $z'$ we have
\begin{align} \label{A-irr}
 w_n(z')&= e^  {(A+Bi)(\ {  {\rm  Log}}\, |z'|+ i(\theta'+ 2n\pi))} 
 \nonumber \\
   &= e^{A( {\rm Log}\ |z'| +i \theta')} \ 
e^{[-B(\theta'+2n \pi)]}e^{i(A2n \pi+ B {\rm Log} |z'|)},  \quad  n\in \mathbb{ Z}.
\end{align}
We first assume $B\ne 0$, e.g., $B>0$. Then
$
 |w_n(z')|= (|z'|^A e^{-B\theta'}) \ e^{-2n B \pi}, \ n\in \mathbb{ Z}$.
 Hence $\lim_{ n\to +\infty} |w_n(z')| = 0$ in $\mathbb{ C}_w$; thus 
$$\lim_{n\to +\infty} (z', w_n(z'))/\sim\  =  [z',0] \in {\bf  T}_a \ \mbox{ in}  \ {\mathcal  
 H}.$$  Since $z'\in \mathbb{ C}^*$ is arbitrary, we have 
$
{  {\bf  T}_{ a}} \subset  \Sigma,
$ 
the closure of $\sigma= \{w= z^{A+Bi}\}/\sim  $ in 
${\mathcal  H}$ .

Since $w=z^{A+Bi}$ can be written as
$$
z= w^{ A'+iB'} \quad \mbox{ where } \ A'={ A}/{(A^2+B^2)}, \ 
B'=-{ B}/{(A^2+B^2)} <0,
$$
we similarly have ${\bf  T}_b \subset \Sigma$. 
This proves 2.(1) in case $B\ne 0$. 

We next assume $B=0$ and $A\ne \rho $. Since the proof is similar, 
we shall prove 2.(1) assuming $-\infty <A< \rho$. For $z\in \mathbb{ C}^*$ we have ${\rm Log} z= \log |z| +i \theta 
\, (0\le \theta<2\pi)$. By analytic continuation of 
$w(z)=z^A=e^{A(\log|z|+i\arg z)}$ along an arbitrary  path 
$l$ from $z$ to $a^kz$ where $k\in \mathbb{ Z}$ is arbitrary, we have 
$$
  w(a^kz)=(a^kz)^A= |a^kz|^Ae^{iA \arg a^kz}= 
|a^kz|^A e^{iA(k \arg a +\theta +2n\pi)}, \quad n\in \mathbb{ Z}.  
$$
Thus $
p_k:= (a^kz, w(a^kz))\in \sigma $. 
In ${\mathcal  H}^*$ the point  $p_k$ coincides with
\begin{align}
 \label{free-k-n}
(z, {  w(a^kz)}/{b^k} )  /\sim &=(z, \widetilde w_k(z))
 /\sim \ \in \sigma 
\end{align}
where $\widetilde w_k(z):=|a^A/b|^k\, e^{ik(A\arg a- \arg b )}e^{i 
A(\theta + 2n\pi )}\in \mathbb{ C}^*_z$.

Using $\rho= \frac{\log |b|}{\log |a|}$,
 \begin{align} \label{free-k}
|\widetilde w_k(z)|= |z|^A (|a|^{kA}/|b|^k)= |z|^A (|a|^{A-\rho})^k.
\end{align}
Since $A<\rho$ and $|a|>1$, it follows that $\lim_{ k \to +\infty}|\widetilde w_k(z)|=0$, so 
that $[z, 0]\in \Sigma $. 
Since $z\in \mathbb{ C}^*$ is arbitrary, we have $\Sigma
\supset {\bf  T}_a$.

Since $w=z^A$ can be written as $z=w^{1/A}$, we have by analytic 
continuation $
q_k:= ((b^k\,w)^{1/A}, b^kw) \in \sigma$ for any $k\in \mathbb{ Z}$.
In ${\mathcal  H}^*$, the point $q_k$ coincides with 
$ ((b^k\,w)^{1/A}/a^k, w)/\sim=:(\widetilde z_k(w), w)/\sim $. 
 Since $|\widetilde z_k(w)|=|w|^{1/A} (|a|^{\rho-A})^{k/A} $, we have 
 $\lim_{ k\to -\infty}|\widetilde z_k(w)| =0$ if $A>0$ and 
$\lim_{ k\to +\infty}|\widetilde z_k(w)| =0$ if $A<0$. Since $w\in 
\mathbb{ C}^*$ is arbitrary, we have $\Sigma  \supset {\bf   T}_b$, which proves 2.(1).

Finally, to prove 2.(2), let $X=\alpha z \frac{\partial }{\partial z} \ne 0$. 
Then the integral curve $ \sigma$  of $X$ passing through $[1,1]$ in 
${\mathcal  H}$ is given by $\{(e^{\alpha t}, 1): t\in \mathbb{ C}\}/\sim \ = \mathbb{ 
C}^*_z \times \{1\}/\sim$. In the fundamental domain 
${\mathcal  F}$, 
$$
 \sigma = (\{0<|z|\le |a|\}, 1) \cup  (\{1< |z|\le |a|\}, 1/b) \cup   (\{1< 
|z|\le |a|\}, 1/b^2) + \ldots,   
$$ 
so that 
$
\Sigma = (\{|z|\le 1\},1)\cup _{n=1}^\infty (\{1\le |z|\le |a|\}, 
1/b^n)\  \cup {\bf  T}_a,
$
proving 2.(2). \hfill $\Box$ 

\vspace{3mm}

 We end this appendix with a remark. Let  $X= \alpha z \frac{\partial 
 }{\partial z}+\beta w \frac{\partial }{\partial w} \not \in \{cX_u: 
 c\in \mathbb{ C}\}  $ with $\alpha \ne 0, \beta \ne 0$ and set 
 $\beta/\alpha= A+Bi$ as in the proof of 2.(1). Fix $(z_0,w_0) \in \mathbb{ C}^* \times \mathbb{ C}^*$ 
 and for $c=w_0/z_0^{\rho}$ consider the integral curve 
$\sigma _c=  \{ w=c z^ {A+Bi}\}/{  \sim} $ of $X$ 
passing through $[z_0,w_0]$ in $ {\mathcal  H}$. For each $z'\in \{1<|z|<|a|\}$ we consider the set  of all points $w_k(z'), \ k=1,2,\ldots  $ 
in $\{|w|<|b|\}$ with $[z',w_k(z')]=(z',w_k(z') \in \sigma _c$. The following fact was used to prove (2) (iii) in Lemma \ref{essen-lemm}.
\begin{remark} \label{al-be-nonzero}
 {\rm  If $A$ is irrational, then  there 
 exists a subsequence $\{w_{k_j}(z')\}_{j=1,2,\ldots} $ 
with the properties that $ \lim_{ j\to \infty}|w_{k_j}(z')|=0$ and 
the closure of the set $\{\arg w_{k_j}(z')\}_{j=1,2,\ldots }$ 
modulo $2\pi $ is equal to $[0, 2\pi]$.}
\end{remark}
\noindent {\bf Proof.} \ Since $\sigma _c= \{w=c 
z^{A+Bi}\}/\sim $ and $\sigma=\{w= z^{A+Bi}\}/\sim$ where 
$\sigma$ is 
defined in the proof of 2.(1), it suffices to prove the result using 
$\sigma _c= \sigma $.  If $B\ne 0$, we can assume $B>0$. Since $A$ is irrational, formula (\ref{A-irr}) gives the result. If $B=0$ we have 
 $A\ne \rho$, and we can assume $-\infty<A<\rho$. In this case, since $A$ is irrational,  formulas (\ref{free-k-n}) and 
 (\ref{free-k}) imply the result.
 \hfill $\Box$

\section{Appendix B: Proof of Lemma \ref{ps-lemma}}

We give the proof of Lemma \ref{ps-lemma}. The lemma is local, hence we may assume from (i) 
and (ii) that the unit outer normal vector of the curve 
$\partial D(0)$ in $\Delta _2$ is $(0,1)$; i.e., $\partial D(0)$ is tangent to 
the $u$-axis at $w=0$ where $w=u+iv$. Thus, we may assume that $\psi(z,w)$ has the 
following Taylor expansion about the origin $(z,w)=(z,(u,v)) =(0,(0,0))$:
\begin{align}
 \label{psi-taylor}
\psi(z,w)&= v+ p_0(z)+ p_1(z)u + p_2(z)u^2+ \ldots =0 \quad 
\end{align}
where each $p_i(z), \ i=0,1,2, \ldots $ is a $C^\omega-$smooth real-valued function and  
$$
p_0(0)=0 \quad \mbox{ and } \quad p_1(0)=0.
$$
We may further assume that formula (\ref{psi-taylor}) holds 
on $(z,u) \in \Delta _1 \times (-r_2,r_2)$ where $\Delta _2=\{|w|<r_2\}$. Thus we write 
\begin{align*}
 D&=\{ v+ p_0(z)+ p_1(z)u + p_2(z)u^2+ \ldots <0: (z,w)\in \Delta_1 \times 
\Delta_2\};\\[3mm]
{\mathcal  S}=\partial D& = 
\{ v+ p_0(z)+ p_1(z)u + p_2(z)u^2+ \ldots =0: (z,w)\in \Delta_1 \times 
\Delta_2\},
\end{align*}
or equivalently, 
\begin{align}\label{D-rep}
D&: \quad  v< -\bigl(p_0(z)+ p_1(z)u + p_2(z)u^2+ \ldots \bigr) \ \ 
 \mbox{ in $\Delta _1 \times \Delta _2$},
\end{align}
and, for each $z\in \Delta _1$, 
$$
S(z): \ \ v=-\bigl(p_0(z)+ p_1(z)u + p_2(z)u^2+ \ldots \bigr)\ \ \ \mbox{ in 
$\Delta _2$}.
$$
In particular, $-ip_0(z)\in S(z)$.  By condition ${(iii)}$ we have 
\begin{align}
 \label{con(iii)}
p_0(z) \not\equiv 0 \qquad \mbox{ on $\Delta _1$.}
\end{align}
Since $\psi(z,w)$ satisfies the Levi condition (\ref{levi}) on $\psi(z,w)=0$,
using the notation 
$$
 \psi(z,w)= \frac{ w-\overline{ w}}{2i}+ p_0(z) +p_1(z)\frac{ w+\overline{ w}}{2}
+p_2(z)(\frac{ w+\overline{ w}}{2})^2+ \ldots,
$$
we calculate to obtain 
\begin{align*}
 {\mathcal  L}\psi(z,w) &=
\left(\frac{\partial^2 p_0(z)}{\partial z \partial \overline{ z}}+ 
\frac{\partial^2 p_1(z)}
{\partial z \partial \overline{ z}}u+
\frac{ \partial p_2(z)}
 {\partial z \partial 
 \overline{ z}} u^2 +\ldots \right)\  
\left| \frac{ 1}{2i}+ \frac{ 1}2  p_1(z) + 
p_2(z) u + \ldots  \right|^2\\
&
- 2 \Re\, \biggl\{  
\left(\frac{ 1}2\frac{\partial p_1(z)}{\partial { z}}+ \frac{\partial 
 p_2(z)}{ \partial \overline{ z} }u + \ldots \right)
\left( \frac{\partial p_0(z)}{\partial \overline{ z}}+ 
\frac{\partial p_1(z)}{\partial \overline{ z}}u
 + \frac{\partial p_2(z)}{\partial \overline{ z}} u^2+ \ldots\right)\\
& \times \left( 
\frac{ 1}{2i}+ \frac{ 1}2  p_1(z) + 
p_2(z) u + \ldots  \right) \biggr\}\\
&
+ \left(\frac{1}2 p_2(z)  + 3p_3(z) u + \ldots 
\right) \left| 
 \frac{\partial p_0(z)}{\partial z}+ \frac{\partial p_1(z)}{\partial z}u
 + \frac{\partial p_2(z)}{\partial z} u^2+ \ldots  \right|^2 \ge 0
\end{align*}
\hspace{7cm}  on $\psi(z,u+iv)=0.$

\noindent In particular,  
\begin{align} 
& \ \  {\mathcal  L}{\psi}(z,0+iv) \nonumber \\
&=
\frac{ 1}4 (1+ p_1(z)^2) \ \frac{\partial^2 p_0(z)}{\partial z 
\partial \overline{ z}}  \nonumber \\ 
& -\frac{ 1}{2} \Re\, \biggl\{
\frac{\partial p_1(z)}{\partial { z}}\ 
\frac{\partial p_0(z)}{\partial \overline{ z}}\ 
(-i+ p_1(z)) \biggr\}+ 
\frac{1}2\  p_2(z) \left|  \frac{\partial p_0(z)}{\partial z}
 \right|^2 \ge 0 \nonumber \\[3mm]
& \hspace{6.5cm} \mbox{on \  $ v+p_0(z)=0$ for $z\in \Delta _1.$}\nonumber  
\end{align}
Since this expression for ${\mathcal  L}{\psi}(z,0+iv)$ is independent of $v$, we have 
\begin{align}\label{levi-2}
(1+ p_1(z)^2) \ \frac{\partial^2 p_0(z)}{\partial z 
\partial \overline{ z}}   -2 & \Re\, \biggl\{ 
\frac{\partial p_1(z)}{\partial { z}} \ 
\frac{\partial p_0(z)}{\partial \overline{ z}}\ 
(-i+ p_1(z)) \biggr\} \nonumber \\ 
& + 
2  p_2(z) \left|  \frac{\partial p_0(z)}{\partial z}
 \right|^2 \ge 0  \qquad \mbox{\quad  for $z \in \Delta _1$.} 
\end{align}
This formula will be used later on in the proof.

\vspace{3mm} 
\noindent  {\it Claim}: To prove the lemma, it suffices to show that for $r_1>0$ sufficiently small  and $\delta _1=\{|z|<r_1\}$, 
 
\begin{center}
 $(\diamondsuit)$ \ \ \ 
{\it  there exists} \ $z^* \in \delta _1$ \ \ {\it   such that} \  \ $p_0(z^*)>0.$
\end{center}

\vspace{3mm} 

\noindent Indeed, if $(\diamondsuit)$ is true, consider the segment $[0,z^*]$ in 
$\delta _1$ and 
the set  
$$ 
{\mathbf  s}:= \bigcup_{z\in [0, z^*]} S(z) \subset \Delta 
_2.
$$
The arc $S(z)$ in $\Delta _2$ varies continuously with $z\in \Delta _1$. Hence it follows from 
$0\in S(0)$,  $-ip_0(z^*) \in S(z^*)$, $-p(z^*) <0$ and (\ref{D-rep}) that 
there exists a sufficiently small disk $\delta _2\subset \Delta _2$ centered at $w=0$ with  $D(0)\cap \delta _2\subset {\mathbf s} $.

\vspace{3mm} 
Thus we turn to the proof of $(\diamondsuit)$. We have two cases, depending on whether $
\frac{\partial p_0}{\partial z}(0)$ vanishes:

\vspace{2mm} 
\noindent {\it Case (i).} \ \  $\displaystyle{ 
\frac{\partial p_0}{\partial z}(0)\ne 0} $.

\vspace{2mm} 
Since $ p_0(0)=0$, we have 
$$
p_0(x,y)= ax+by+ O(|z|^2) \quad \mbox{ near $z=0$} 
$$ 
with $(a,b)\ne (0,0)$. It is clear that there exist $z^*\in \delta _1$ 
 which satisfies $(\diamondsuit)$. 

\vspace{2mm} 
\noindent {\it Case (ii).} \ \ $ \displaystyle{ \frac{\partial p_0}{\partial z}(0)=0} $.

\vspace{2mm} 

\noindent In this case,
we have the following Taylor expansion of $p_0(z)$ about $z=0$:
\begin{align*} 
 (1) \quad   p_0(z) &= 
 \Re\, \{a_{20} z^2\} + a_{11}z \overline{ z} \\[2mm]
& \quad + 
\ldots + \Re\,  J_{2n-1} + \Re\, J_{2n} + O(|z|^{2n+1}) \ \ \mbox{ near $z=0$},
\end{align*}
where 
$$
J_{2n-1}= \Re\, \bigl\{\sum_ {k=0}^n  a_{2n-k,k}\,z^{2n-k} \overline{ z}^k\bigr\}, \quad  
J_{2n}= \Re\, \bigl\{  \sum_ {k=0}^n  a_{2n-k,k}\,z^{2n-k}\overline{ z}^k
\bigr\}+ a_{nn}|z|^{2n}.
$$
Here $a_{ij}$ is, in general, a complex number for $i\ne j$; while $a_{ii}$ is real.

\vspace{2mm}  
\noindent {\it  $1^{st}$ step}: Since $\frac{\partial 
p_0}{\partial z}(0)=0$ and $p_0(0)=p_1(0)=0$,  inequality (\ref{levi-2}) 
reduces to 
$$
\frac{\partial^2 p_0}{\partial z \partial \overline{ z}  }(0) \ge 0, 
\quad \mbox{ i.e.,} \quad a_{11}\ge 0.
$$

If $a_{11}>0$, (1) implies that
\begin{align*}
 \frac{\partial^2 p_0}{\partial z \partial \overline{ z}  }(z)&=
a_{11} +O(|z|)\ge \frac{ a_{11}}2 >0 \quad \mbox{ near $z=0$}.
\end{align*}
Thus $p_0(z)$ is strictly subharmonic on a sufficiently small 
disk $\delta _1':= \{|z|<r'\} \subset \delta _1$; hence there exists $z^*$ with 
$|z^*|= \frac{ r'}2$ and $p_0(z^*)>p_0(0)=0$, proving $(\diamondsuit)$.

If $a_{11}=0$, then (1) becomes, for $z=r e^{i\theta}$, 
\begin{align*}
  p_0(z) &=
 \Re\, \{a_{20} z^2\}  +O(|z|^3)=
|z|^2 \ \Re\, \bigl\{a_{20}e^{2i\theta} +O(|z|) \bigr\} 
\ \  \mbox{ near $z=0$}.
\end{align*}

If $a_{20}=|a_{20}| \, e^{i\theta_0}\ne 0$, then  
for $z^*\in \delta _1$ of the form $z^*=
r^*e^{-i \theta_0/2} \ne 0$ with $r^*$ sufficiently small, we have 
$$
 p_0(z^*)= (r^*)^2 \ \bigl( \ |a_{20}| +O(|z^*|)\bigr) \ge (r^*)^2\ 
\frac{ |a_{20}|}2 >0,
$$
which proves ($\diamondsuit$).

\vspace{2mm} 
Thus it suffices to prove $(\diamondsuit)$ in the following two 
cases when $n\ge 2$:
$${\rm Case\ (I)}: \  p_0(z)= J_{2n-1}(z) +O(|z|^{2n}) \ \mbox{  near $z=0$}
$$
where
$$
J_{2n-1}(z):= \Re\, \{ a_{2n-1}z^{2n-1}+ a_{2n-2}
z^{2n-2} \overline{ z} + \ldots + a_{n}z^n \overline{ z}^{n-1} \} 
\ \ \mbox{in }\ \mathbb{ C}_z;
$$
$a_{i}$ is, in general, a complex number; and
\begin{align}
 \label{non-01}
 (a_{2n-1}, a_{2n-2}, \ldots , a_{n}) \ne (0,0, \ldots ,0).
\end{align}
$${\rm Case\ (II)}: \quad  p_0(z)= J_{2n}(z) +O(|z|^{2n+1}) \ \mbox{  near $z=0$}
$$
where
$$
J_{2n}(z):= \Re\, \{ a_{2n}z^{2n-1}+ a_{2n-1}
z^{2n-1} \overline{ z} + \ldots + a_{n+1}z^{n+1} \overline{ z}^{n-1} \} 
+a_n |z|^{2n}
\ \ \mbox{in }\ \mathbb{ C}_z;
$$
$a_{i}$ for $n+1\le i\le 2n$ is, in general, a complex number; $a_{n}$ is a 
real number; and
\begin{align}
 \label{non-02}
 (a_{2n}, a_{2n-1}, \ldots , a_{n+1}, a_n) \ne (0,0, \ldots ,0,0).
\end{align}

We first assume ${\rm Case \ (I)}$. Setting $z=|z|e^{i\theta}$, we 
 have 
\begin{align*}
 J_{2n-1}(z)&=|z|^{2n-1} 
\ \Re\, \{ a_{2n-1}e^{i(2n-1)\theta}+ a_{2n-2}
e^{i(2n-3)\theta}+ \ldots + a_{n}e^{i\theta}\}  
\ \ \mbox{in }\ \mathbb{ C}_z
\end{align*}
We consider the polynomial in $Z$ defined by 
$$
g(Z):= a_{2n-1} Z^{2n-1} + a_{2n-2}Z^{2n-3}+ \ldots + a_{n}Z.
$$
Note that $g(Z)\not \equiv 0$ by (\ref{non-01}). Thus $g(Z)\ne 0$ for all 
$Z$ with $|Z|=r$ for some $0<r<1$. Since $g(0)=0$, by the 
argument principle $\int_{|Z|=r }d \arg\, g(Z)\ge 1$, hence there exists $0\le \theta'< 2\pi $ such that  $\Re\, 
g(re^{i\phi'})>0$. 
By the maximum principle for the harmonic function $\Re\, g(Z)$ on 
$\{|Z|\le 1\}$, there exists $0\le \theta^*\le 2\pi$ such that 
 $$
A:= \Re\, g(e^{i\theta^*}) \ge \Re\, g(re^{i\theta'})>0.
$$
Since $J_{2n-1}(z)= |z|^{2n-1}g(e^{i\theta})$, we have 
\begin{align*}
 p_0(|z|e^{i\theta^*})&= |z|^{2n-1} A +O(|z|^{2n}) \quad \mbox{ for $0<|z| 
\ll 1$}\\
&\ge  |z|^{2n-1} A/2 >0 \quad \mbox{ for $0<|z|\ll 1$},
\end{align*}
showing that $(\diamondsuit)$ is true in ${\rm Case \ (I)}$. 

We next assume ${\rm Case \ (II)}$. For $z=|z|e^{i\theta}$
\begin{align}\label{umai2}
 \frac{\partial^2 p_0(z)}{\partial z \partial \overline{ z}}
&= |z|^{2n-2} \biggl[ \Re\, \{(*)\} \ +  a_{n} + O(|z|)\biggr]
\end{align}
where 
\hspace{ -1cm}
$$
(*)=
(2n-1)a_{2n-1}e^{i(2n-2)\theta}+(2n-2)2\cdot a_{2n-2}e^{i(2n-4)\theta}
+ \ldots + a_{n+1}e^{i2\theta}. 
$$
We substitute this in (\ref{levi-2}) to obtain 
\begin{align*}
 &\bigl(1+O(1)^2\bigr) |z|^{2n-2}
\biggl(
\biggl[ \Re\,  \{(*) \}  + n^2a_n \biggr]
+O(|z|)\biggr)\\
&
-2 \Re\, \biggl\{ O(1) O(|z|^{2n-1}) (-i+O(1)) \biggr\}\ 
 + 2\,O(|z|)\, O(|z|^{2n-1})^2 
\ge 0 
\end{align*}
for $|z|$ sufficiently small. Dividing both sides by $(1+ O(1)^2)\,|z|^{2n-2}>0$ with $|z|>0$ and then 
letting $|z|\to 0$, we have 
\begin{align}
 \label{umai}
\Re\, \{ (*) \} + n^2\,a_n \ge 0 \quad \mbox{ for all} \ 0\le 
 \theta< 2\pi.
\end{align} 
We substitute this in the definition of $p_0(z)$ in Case (II) to obtain
\begin{align*}
&p_0(z)\ge |z|^{2n}\  \Re\,
\biggl\{ 
a_{2n} e^{i2n\theta} 
+ a_{2n-1}
\bigl(1-\frac{2n-1}{n^2} \bigr) e^{i(2n-2)\theta}
+\\
& a_{2n-2}
\bigl(1-\frac{(2n-2)2}{n^2}\bigr) e^{i(2n-4)\theta}
+ \ldots +
a_{n+1}\bigl(
1-\frac{ (n+1)(n-1)}{n^2}\bigr) 
e^{i2\theta} 
\biggr\}\
 \\
 &+  O(|z|^{2n+1}) 
\end{align*}
for $|z|$ sufficiently small. 

We divide the proof of Case (II) in two subcases:
\begin{align*}
 &\mbox{ Case (II-1):} \quad  (a_{2n}, a_{2n-1}, \ldots , a_{n+1}) \ne (0,0, 
\ldots , 0);\\
&\mbox{ Case (II-2):} \quad  (a_{2n}, a_{2n-1}, \ldots , a_{n+1}) =(0,0, 
\ldots , 0).
\end{align*}
From (\ref{non-02}), $a_{n}\ne 0$ in Case (II-2). In Case (II-1) we consider the polynomial 
\begin{align*}
 g(Z)&= 
a_{2n} Z^{2n} 
+ a_{2n-1}
\bigl(1-\frac{2n-1}{n^2} \bigr) Z^{2n-2}
+
a_{2n-2}
\bigl(1-\frac{(2n-2)2}{n^2}\bigr) Z^{2n-4}
 \\
& \qquad \ldots +
a_{n+1}\bigl(
1-\frac{ (n+1)(n-1)}{n^2}\bigr) 
Z^{2}.
\end{align*}
Since $n\ge 2$, we have $(1- \frac{ (2n-k)k}{n^2}) \ne 0$ for $k=1,2,\ldots , 
n-1$ so that $g(Z) \not \equiv 0$ on $\mathbb{ C}_Z$ and $g(0)=0$. 
By the same reasoning as in Case (I) we have the existence of $0\le 
\theta^*<2\pi$ and $A>0$ with  
$$
p_0(|z|e^{i\theta^*}) \ge |z|^{2n} \, A/2 >0 \quad \mbox{ for 
$0< |z|\ll 1$},
$$
which proves $(\diamondsuit)$ in Case (II-1). 

In Case (II-2) we have $(*)=0$ in (\ref{umai2}) and hence $a_n\ge 0$ 
from (\ref{umai}); thus $a_n>0$. Using (\ref{umai2}) we have 
$$
\frac{\partial^2 p_0(z) }{\partial z \partial \overline{ z}}
\ge 
|z|^{2n-2} a_n +O(|z|^{2n-2}) \ge |z|^{2n-2} a_n/2 \ge 0
$$
for $z$ in a sufficiently small disk $\delta $ centered at $z=0$. In other words, $p_0(z)$ is subharmonic on $\delta $ and is strictly subharmonic 
in $\delta \setminus \{0\}$. Thus, for a given  $0<r<r_0$, we can find $0\le 
\theta^*<2\pi$ with $p_0( r e^{i\theta^*})>0$, which proves 
$(\diamondsuit)$ in Case (II-2). This completes the proof of $(\diamondsuit)$.
 \hfill $\Box$

\end{document}

lllllllllllllllllllllllllllllllllllll

 Princeton Mathematical Series, No. 26, Princeton Univ. Press,
  Princeton, 1960.}

\end{thebibliography}
\end{document}